\documentclass[11pt]{article}
\usepackage{amssymb}
\usepackage{amsmath}
\usepackage{amsthm}
\newcommand{\cA}{{\cal A}}

\newcommand{\cD}{{\cal D}}
\newcommand{\cH}{{\cal H}}

\newcommand{\cO}{{\cal O}}
\newcommand{\cL}{{\cal L}}

\newcommand{\cF}{{\cal F}}

\newcommand{\cP}{{\cal P}}

\newcommand{\cU}{{\cal U}}
\newcommand{\cV}{{\cal V}}
\newcommand{\cW}{{\cal W}}

\newcommand{\cY}{{\cal Y}}

\renewcommand{\AA}{{\mathbb A}}
\newcommand{\GG}{{\mathbb G}}

\newcommand{\ZZ}{{\mathbb Z}}
\newcommand{\QQ}{{\mathbb Q}}

\renewcommand{\gg}{\mathfrak{g}}

\newcommand{\gp}{\mathfrak{p}}
\newcommand{\gq}{\mathfrak{q}}

\newcommand{\gA}{\mathfrak{A}}

\newcommand{\on}{\operatorname}

\newcommand{\Rep}{{\on{Rep}}}

\newcommand{\Qlb}{\mathbb{\bar Q}_\ell}
\newcommand{\Gm}{\mathbb{G}_m}
\newcommand{\Ga}{\mathbb{G}_a}
\newcommand{\A}{\mathbb{A}}

\newcommand{\toup}[1]{\stackrel{#1}{\to}}
\newcommand{\hook}[1]{\stackrel{#1}{\hookrightarrow}}
\newcommand{\df}{\stackrel{\rm{def}}{=}}
\newcommand{\getsup}[1]{\stackrel{#1}{\gets}}
\newcommand{\Sp}{\on{\mathbb{S}p}}
\newcommand{\Spin}{\on{\mathbb{S}pin}}
\newcommand{\GSpin}{\on{G\mathbb{S}pin}}
\newcommand{\GSp}{\on{G\mathbb{S}p}}

\newcommand{\CT}{\on{CT}}
\newcommand{\Hom}{\on{Hom}}

\newcommand{\End}{\on{End}}

\newcommand{\Sym}{\on{Sym}}
\newcommand{\SO}{\on{S\mathbb{O}}}

\newcommand{\RG}{\on{R\Gamma}}

\newcommand{\Spr}{{{\cal S}pr}}
\newcommand{\length}{\on{\ell}}

\newcommand{\Bun}{\on{Bun}}
\newcommand{\Bunb}{\on{\overline{Bun}} }

\newcommand{\Spec}{\on{Spec}}
\newcommand{\Specf}{\on{Spf}}

\newcommand{\supp}{\on{supp}}

\newcommand{\HOM}{{{\cal H}om}}

\newcommand{\Gr}{\on{Gr}}
\newcommand{\Grb}{\overline{\Gr}}

\newcommand{\GL}{\on{GL}}
\newcommand{\PSL}{\on{PSL}}

\newcommand{\id}{\on{id}}

\newcommand{\QED}{$\square$} 
\newcommand{\Fq}{\mathbb{F}_q}  
\newcommand{\iso}{{\widetilde\to}}

\newcommand{\comp}{\circ}

\renewcommand{\H}{{\on{H}}}   
\newcommand{\R}{\on{R}\!}   

\newcommand{\DD}{\mathbb{D}}  
\newcommand{\D}{\on{D}}       
\newcommand{\wt}{\widetilde}

\newcommand{\select}[1]{{\it{#1}}}

\newcommand{\K}{{\on{K}}}

\newcommand{\<}{\langle}
\renewcommand{\>}{\rangle}

\newcommand{\Av}{\on{Av}}
\newcommand{\ev}{\on{ev}}
\newcommand{\conv}{\on{conv}}
\newcommand{\Loc}{\on{Loc}}
\newcommand{\Lie}{\on{Lie}}
\newcommand{\Sph}{\on{Sph}}
\newcommand{\Res}{\on{Res}}
\newcommand{\leM}{{\atop{\stackrel{\le}{M}}}}
\newcommand{\ttimes}{\tilde\times}
\newcommand{\Xo}{\breve{X}}

\newtheorem{Lm}{Lemma}

\newtheorem{Pp}{Proposition}
\newtheorem{Cor}{Corolary}
\newtheorem{Con}{Conjecture}
\newtheorem{Slm}{Sublemma}

\newcommand{\I}{\cal I} 

\theoremstyle{remark}
\newtheorem{Rem}{Remark}
\newtheorem{Rems}{Remarks}

\theoremstyle{definition}
\newtheorem{Def}{Definition}

\newenvironment{Prf}{\par\noindent {\it Proof }}{\QED}

\begin{document}

\author{Sergey Lysenko}
\title{On automorphic sheaves on $\Bun_G$}
\date{}
\maketitle
\begin{abstract}
\noindent{\scshape Abstract}\hskip 0.8 em 
Let $X$ be a smooth projective connected curve over an algebraically closed field
$k$ of positive characteristic. Let $G$ be a reductive group over $k$,
$\gamma$ be a dominant coweight for $G$, and $E$ be an $\ell$-adic $\check{G}$-local system on
$X$, where $\check{G}$ denotes the Langlands dual group (over
$\Qlb$). Let $\Bun_G$ be the moduli stack of $G$-bundles on $X$.

 Under some conditions on the triple $(G,\gamma,E)$ we propose 
a conjectural construction of a distinguished 
$E$-Hecke automorphic sheaf on $\Bun_G$. 
We are motivated by a construction of automorphic forms suggested by
Ginzburg, Rallis and Soudry in \cite{GRS1, GRS2}.
We also generalize Laumon's theorem (\cite{La2}, Theorem~4.1) for
our setting. Finally, we formulate an analog of the Vanishing Conjecture of Frenkel,
Gaitsgory and Vilonen for Levi subgroups of $G$.
\end{abstract} 

{\centerline{\scshape 1. Introduction}}

\medskip\noindent
 Let $F$ be a number field, $\AA$ be its ring of adeles. If
$G$ is one of the groups $\SO_{2n+1}$, $\Sp_{2n}$ or $\SO_{2n}$ then consider standard
representation
$\check{G}\to \check{H}$ of the Langlands dual group $\check{G}$, so here $H$ is $\GL_{2n}$,
$\GL_{2n+1}$ or $\GL_{2n}$ respectively. For an irreducible, automorphic, cuspidal
representation $\tau$ of $H(\AA)$ satisfying some additional conditions, D. Ginzburg, S. Rallis
and D. Soudry have proposed a conjectural construction of an irreducible, automorphic cuspidal
representation $\sigma$ of $G(\AA)$ which lifts to $\tau$ (cf., \cite{GRS1, GRS2}).

 For example, consider $G=\SO_{2n+1}$.  Let $X$ be a smooth projective
absolutely irreducible curve over $\Fq$. Consider the Langlands dual group $\check{G}=\Sp_{2n}$
over $\Qlb$. Let $H=\GL_{2n}$ over $\Fq$, $\check{H}=\GL_{2n}$ over $\Qlb$. The standard
representation $V$ of $\check{G}$ is a map $\check{G}\to\check{H}=\GL(V)$. 

  Let $E$ be an $\ell$-adic $\check{G}$-local system on $X$, assume that $V_E$ is irreducible.
According to \cite{Laf}, theorem VII.6, irreducibility of $V_E$ implies that $\End(V_E)$ is pure
of weight zero. It follows that for each closed point $x\in X$ the local L-function $L(E_x,
\check{\gg}, s)$ is regular at $s=1$, and the corresponding irreducible unramified
representation of $G(F_x)$ is generic. Here
$F_x$ denotes the completion of $\Fq(X)$ at $x$, and $\check{\gg}=\Lie \check{G}$.
An analog of D.~Ginzburg, S.~Rallis and D.~Soudry's conjecture for function field is to predict
that in the L-packet of automorphic forms corresponding to $E$ there exists a unique
nonramified cuspidal generic form $\varphi_E:\Bun_G(\Fq)\to\Qlb$ 
(cf. \cite{GRS2}, Conjecture on p. 809 and \cite{GRS1}). 

 In \select{loc.cit.} an additional condition is required: the L-function $L(E,\wedge^2 V,s)$
has a pole of order exactly one at $s=1$. This condition is satisfied in our situation. Indeed,
$\wedge^2 V=V'\oplus\Qlb$, where $V'$ is an irreducible representation of
$\check{G}$. Since $V$ is self-dual, $\H^0(X\otimes\bar\Fq, V'_E)=0$
and the L-function $L(E,V',s)$ is a polynomial in $q^{-s}$. The purity argument shows that
$L(E,V',1)\ne 0$.
 
 In this paper we consider the problem of constructing a geometric counterpart of
$\varphi_E$. Given a reductive group $G$, a dominant coweight $\gamma$ and a $\check{G}$-local 
system $E$ on $X$, we impose on these data some conditions similar to the above. Then we
propose a conjectural construction of a distinguished $E$-Hecke eigensheaf on the moduli stack
$\Bun_G$ of $G$-bundles on $X$. Our approach applies to root systems $A_n, B_n, C_n$ for all
$n$, $D_n$ for odd $n$, and also $E_6, E_7$. For $\GL_n$ our method reduces to the
one proposed by Laumon in \cite{La}. 

 The construction is exposed in Sect.~2, 3. In Sect.~4 we study the additional
structure on Levi subgroups induced by $\gamma$, and prove a generalization of Laumon's
theorem (\cite{La2}, Theorem~4.1) for our setting. We discuss its applications
to cuspidality and formulate an analog of the Vanishing Conjecture of Frenkel, Gaitsgory
and Vilonen for Levi subgroups of $G$.

\bigskip
\centerline{\scshape 2. Statements and conjectures}

\bigskip\noindent
2.1 {\scshape Notation\ } Throughout, $k$ will denote an algebraically closed field of
characteristic $p>0$. Let
$X$ be a smooth projective connected curve over $k$. Fix a prime $\ell\ne p$. For a 
$k$-scheme (or $k$-stack) $S$ write $\D(S)$ for the bounded derived category of $\ell$-adic
\'etale sheaves on $S$. 

 Let $G$ be a connected reductive group over $k$. Fix a Borel subgroup $B\subset G$. Let $N\subset B$ be its unipotent radical and $T=B/N$ be the
"abstract" Cartan. Let $\Lambda$ denote the coweight lattice. The weight lattice is denoted by
$\check{\Lambda}$. The semigroup of dominant coweights (resp., weights) is denoted $\Lambda^+$
(resp., $\check{\Lambda}^+$). The set of vertices of the Dynkin diagram of $G$ is denoted by
$\I$. To each $i\in\I$ there corresponds a simple root $\check{\alpha}_i$ and a simple coroot
$\alpha_i$. By $\check{\rho}\in\check{\Lambda}$ is denoted the half sum of positive roots of
$G$ and by $w_0$ the longest element of the Weil group $W$. 

 Let $\Lambda^{pos}$ denote
$\ZZ_+$-span of positive coroots. 
The set $\Lambda^+$ is equiped with the order $\nu_1\le
\nu_2$ iff $\nu_2-\nu_1\in \Lambda^{pos}$. Similarly, we have an order on
$\check{\Lambda}^+$.
 
 To a dominant weight $\check{\lambda}$ one attaches the Weil $G$-module
$\cV^{\check{\lambda}}$ with a fixed highest weight vector $v^{\check{\lambda}}\in
\cV^{\check{\lambda}}$. For any pair $\check{\lambda},\check{\nu}\in\check{\Lambda}^+$ there is
a canonical map
$\cV^{\check{\lambda}+\check{\nu}}\to\cV^{\check{\lambda}}\otimes\cV^{\check{\nu}}$ sending
$v^{\check{\lambda}+\check{\nu}}$ to $v^{\check{\lambda}}\otimes v^{\check{\nu}}$. 

 For $\lambda\in\Lambda^+$ write $V^{\lambda}$ for the irreducible representation of $\check{G}$ of highest weight $\lambda$.

  The trivial $G$-bundle on a scheme is denoted by $\cF^0_G$. Recall that for any finite
subfield $k'\subset k$ and any non-trivial character $\psi: k'\to\Qlb$ one can construct the
Artin-Shrier sheaf $\cL_{\psi}$ on $\GG_{a,k}$. The intersection cohomology
sheaves are normalized to be pure of weight zero.

\bigskip\noindent
2.2 {\scshape Additional data and assumptions\ } We say that $\gamma\in\Lambda^+$ is \select{minuscule} if $\gamma$ is a minimal element of $\Lambda^+$ and $\gamma\ne 0$. 
If $\gamma\in\Lambda^+$ is minuscule then, by (Lemma~1.1, \cite{NP}), for any root $\check{\alpha}$ we have
$\<\gamma,\check{\alpha}\>\in \{0,\pm 1\}$, and the set of weights of $V^{\gamma}$ coincides with the $W$-orbit of $\gamma$.  For example, if $\gamma\ne 0$ is orthogonal to all roots then $\gamma$ is minuscule. 
One checks that the natural map from the set of minuscule dominant coweights to $\pi_1(G)$ is injective.

\begin{Def} We say that $\{\gamma\}$ is a \select{1-admissible datum} if the following
conditions hold
\begin{itemize}
\item the center $Z(G)$ is a connected 1-dimensional torus;
\item $\pi_1(G)\iso\ZZ$;
\item $\gamma\in\Lambda^+$ is a minuscule dominant coweight whose image $\theta$ in
$\pi_1(G)$ generates $\pi_1(G)$;
\item $V^{\gamma}$ is a faithful representation of $\check{G}$.
\end{itemize}
\end{Def}

Fix a 1-admissible datum $\gamma$. 
For $k\ge 0$ set $\Lambda^{+,k\theta}_{G,S}=\{\mu\in\Lambda^+\mid \; \mu\le
k\gamma\}$. Let $\Lambda^+_{G,S}$ be the union of $\Lambda^{+,k\theta}_{G,S}$ for
$k\ge 0$. Set
$$
\check{\Lambda}^+_S=\{\check{\lambda}\in\check{\Lambda}^+\mid \;
\<w_0(\lambda),\check{\lambda}\>\ge 0 \;\mbox{for any}\;
\lambda\in\Lambda^+_{G,S}\}
$$
Let $\check{\omega}_0$ be the generator of the group
of weights orthogonal to all coroots, we fix $\check{\omega}_0$ by requiring $\<\theta,
\check{\omega}_0\>=1$. For $i\in\cal I$ denote by $\check{\omega}_i\in\check{\Lambda}^+$ the
fundamental weight corresponding to $\alpha_i$ that satisfies
$\<w_0(\gamma),\check{\omega}_i\>=0$. Note that $\check{\omega}_0, \check{\omega}_i (i\in\cal
I)$ form a basis of $\check{\Lambda}$.

 The following lemma is straightforward. 

\begin{Lm}
The semigroup $\check{\Lambda}^+_S$ is the $\ZZ_+$-span of $\check{\omega}_0,
\check{\omega}_i (i\in\cal I)$. Besides
$$
\Lambda^+_{G,S}=\{\lambda\in\Lambda^+\mid \; \<w_0(\lambda),\check{\lambda}\>\ge 0 \;
\mbox{for all}\; \check{\lambda}\in \check{\Lambda}^+_S\} \eqno{\square}
$$
\end{Lm}

Since $\pi_1(G)\iso\ZZ$, it follows that $[G,G]$ is simply-connected.
Note that for $\check{\mu}\in\check{\Lambda}^+, \check{\lambda}\in\check{\Lambda}^+_S$
the condition $\check{\mu}\le\check{\lambda}$ implies 
$\check{\mu}\in\check{\Lambda}^+_S$. Note that $\{-w_0(\gamma)\}$ is also a
1-admissible datum.

 Since $V^{\gamma}$ is faithful, the weights of $V^{\gamma}$ generate $\Lambda$ and for each
$i\in\cal I$ we have $\<\gamma,\check{\omega}_i\> >0$. For each maximal positive root
$\check{\alpha}$ we have $\<\gamma,\check{\alpha}\>=1$. In particular, if the root system of $G$ is irreducible (so, nonempty) then $\gamma$ is a fundamental coweight corresponding to some simple root.

 Some examples of 1-admissible data are given in the appendix.

\bigskip\noindent
2.2.1 Consider the formal disk $\cD=\Specf(k[[t]])$. Recall that the Affine Grassmanian $\Gr_G$ is
the ind-scheme classifying pairs $(\cF_G,\beta)$, where $\cF_G$ is a $G$-bundle on $\cD$ and
$\beta:
\cF_G\iso\cF^0_G$ is a trivialization over the punctured disk $\cD^*=\Spec k((t))$. 
Define the positive part $\Gr_G^+\subset\Gr_G$ of $\Gr_G$ as a closed
subscheme given by the following condition:

  $\cF_G\in\Gr_G^+$ if for every $\check{\lambda}\in\check{\Lambda}^+_S$ the map
$$
\beta^{\check{\lambda}}: \cV^{\check{\lambda}}_{\cF_G}\mid_{\cD^*}\to 
\cV^{\check{\lambda}}_{\cF^0_G}\mid_{\cD^*}
$$
is regular on $\cD$. Note that $\Gr_G^+$ is invariant under the natural action of $G(k[[t]])$.

Recall that for $\mu\in\Lambda^+$ one has the closed subscheme
$\Grb^{\mu}_G\subset\Gr_G$ (cf.
\cite{BG}, sect. 3.2). One checks that 
$\Grb^{\mu}_G\subset \Gr_G^+$ iff $\mu\in\Lambda^+_{G,S}$. Let $\pi_1^+(G)\subset \pi_1(G)$
be the image of $\Lambda_{G,S}^+$ under the projection $\Lambda\to \pi_1(G)$. 

 For $\nu\in\pi_1(G)$ the connected component $\Gr_G^{\nu}$ of
$\Gr_G$ is given by the condition:
$$
\cV^{\check{\omega}_0}_{\cF^0_G}(-\<\nu,\check{\omega}_0\>)\,\iso\,
\cV^{\check{\omega}_0}_{\cF_G}
$$
For $\nu\in\pi_1(G)$ set $\Gr_G^{+,\nu}=\Gr_G^+\cap \Gr_G^{\nu}$.

\bigskip\noindent
2.3 Denote by $\Bun_G$ the moduli stack of $G$-bundles on $X$. 
Let $\cH_G^+$ be the corresponding positive part of the Hecke stack, it classifies
collections: $\cF_G,\cF'_G\in\Bun_G$, an effective divisor $D$ on $X$, an isomorphism 
$\beta: \cF_G\mid_{X-D}\iso\cF'_G\mid_{X-D}$ such that for each
$\check{\lambda}\in\check{\Lambda}^+_S$ the map 
$$
\beta^{\check{\lambda}}: \cV^{\check{\lambda}}_{\cF_G}\hook{}\cV^{\check{\lambda}}_{\cF'_G}
$$
extends to an inclusion of coherent sheaves on $X$, 
and $\cV^{\check{\omega}_0}_{\cF_G}(D)\,\iso\,
\cV^{\check{\omega}_0}_{\cF'_G}$. 
For $k\ge 0$ let $\cH^{+,k}_G\subset \cH^+_G$ be given by
$\deg D=k$. 
 
\bigskip\noindent
2.4 {\scshape Version of Laumon's sheaf\ }  Given a local system $W$ on $X$ and $d\ge 0$, one
defines a sheaf $\cL^d_W$ on $\cH^{+,d}_G$ as follows. 

\smallskip

 Let $\tilde\cH^{+,d}_G$ be the stack
of collections: $(\cF^1_G,\ldots,\cF^{d+1}_G)$ and $(\beta^i, x_i\in X)_{i=1,\ldots, d}$,
where $\cF^i_G\in\Bun_G$ and
$$
\beta^i: \cF^i_G\mid_{X-x_i}\,\iso\,\cF^{i+1}_G\mid_{X-x_i}
$$
is an isomorphism such that $(\cF^i_G,\cF^{i+1}_G,\beta^i, x_i)\in \cH^{+,1}_G$ for
$i=1,\ldots, d$. 

Note that $\tilde\cH^{+,d}_G$ is smooth, because $\gamma$ is minuscule. We have a
convolution diagram 
$$
\begin{array}{ccc}
\tilde\cH^{+,d}_G & \toup{\wt\supp} & X^d\\
\downarrow\lefteqn{\scriptstyle p} && \downarrow\\
\cH^{+,d}_G & \toup{\supp} & X^{(d)}
\end{array}
$$
where $\wt\supp$ (resp., $p$) sends the above collection to $(x_1,\ldots,x_d)$ (resp., to
$(\cF^1_G,\cF^{d+1}_G,\beta,D)$ with $D=x_1+\ldots+x_d$.  Let $^{rss}X^{(d)}\subset
X^{(d)}$ be the open subscheme classifying reduced divisors. Over $^{rss}X^{(d)}$, this
diagram is cartesian. 

 The following proposition is an
immediate corolary of (Lemma 9.3 \cite{NP}).

\begin{Pp} The map $p$ is representable proper surjective and small. \QED
\end{Pp}

 Set $\Spr_W^d=p_!\wt\supp^*W^{\boxtimes\, d}[m](\frac{m}{2})$ with $m=\dim\cH^{+,d}_G$. This
is a perverse sheaf, the Goresky-MacPherson extension from ${\supp}^{-1}({^{rss}X}^{(d)})$. It
is equiped with a canonical action of $S_d$. Define $\cL^d_W$ to be the $S_d$-invariants of
$\Spr^d_W$. We have $\DD(\Spr^d_W)\iso \Spr^d_{W^*}$ canonically and $\DD(\cL^d_W)\iso
\cL^d_{W^*}$.

 We have a diagram 
$$
\Bun_G\getsup{\gp}\cH^+_G\toup{\gq}\Bun_G
$$
where $\gp$ (resp., $\gq$) sends $(\cF_G,\cF'_G,\beta)$ to $\cF_G$ (resp., $\cF'_G$). 
By (property 3, sect. 5.1.2 \cite{BG}), the sheaf $\cL^d_W$ is ULA with respect
to both projections $\gp$ and $\gq$.

 Let $r=\dim V^{\gamma}$. For a partition $\mu=(\mu_1\ge \ldots\ge\mu_r\ge 0)$ of $d$ define
the polynomial functor $W\mapsto W_{\mu}$ of a $\Qlb$-vector space $W$ by
$$
W_{\mu}=(W^{\otimes d}\otimes U_{\mu})^{S_d},
$$
where $U_{\mu}$ stands for the irreducible representation of $S_d$ corresponding to $\mu$. 
For $d>0$ let $\length(\mu)$ be the greatest index $i\le r$ such that $\mu_i\ne 0$. For $d=0$
let $\length(\mu)=0$. If $\length(\mu)$ is less or equal to $\dim W$ then $W_{\mu}$ is the
irreducible representation of $\GL(W)$ with h.w. $\mu$, otherwise it vanishes.

 For $\nu\in\Lambda^+$ let $\cA_{\nu}$ denote the IC-sheaf on $\Grb^{\nu}$. Recall that the
category $\Sph(\Gr_G)$ of spherical perverse sheaves on $\Gr_G$ consists of direct sums of 
$\cA_{\nu}$, as $\nu$ ranges over the set of dominant coweights.
We have the Satake equivalence of tensor categories $\Loc: \Rep(\check{G})\to
\Sph(\Gr_G)$ (cf. Theorem~3.2.8, \cite{BG}). In particular, we have
$\Loc(V^{\nu})=\cA_{\nu}$.   

 Consider a $k$-point $(\sum d_kx_k,\, \cF'_G)$ of $X^{(d)}\times\Bun_G$.
The fibre of $\supp\times\gq:\cH^{+,d}_G\to X^{(d)}\times\Bun_G$ over
this point identifies with 
\begin{equation}
\label{fibre_1}
\prod_k \Grb^{d_k\gamma}_G
\end{equation}
 Set $d_G=\dim\Bun_G$. The Satake equivalence yields the following description. 
\begin{Pp} 
\label{Pp_description_Laumon}
For any local system $W$ on $X$ the restriction of $\cL^d_W$ to the fibre
(\ref{fibre_1}) of $\supp\times\gq$ identifies with the exterior product $(\boxtimes_{k}
{\cL_k})[d+d_G](\frac{d+d_G}{2})$, where each $\cL_k$ is
$$
\cL_k=\mathop{\oplus}\limits_{\mu} \Loc((V^{\gamma})_{\mu})\otimes (W_{x_k})_{\mu},
$$
the sum being taken over the set of partitions of $d_k$ of length $\le r$. \QED
\end{Pp}

\smallskip
\noindent
2.5 Given a local system $W$ on $X$, for $d\ge 0$ define a
functor $\Av^d_W:\D(\Bun_G)\to\D(\Bun_G)$ by
$$
\Av^d_W(K)=\gq_!(\gp^*K\otimes\cL^d_W)[-d_G](\frac{-d_G}{2})
$$
Let also $\Av^{-d}_W:\D(\Bun_G)\to\D(\Bun_G)$ be given by
$$
\Av^{-d}_W(K)=\gp_!(\gq^*K\otimes\cL^d_W)[-d_G](\frac{-d_G}{2})
$$
The functors $\Av^d_W$ and $\Av^{-d}_{W^*}$ are both left and right adjoint to each other.
As in (Proposition 9.5, \cite{FGV2}) one proves

\begin{Pp} 
\label{Pp2}
Let $K$ be a Hecke eigensheaf on $\Bun_G$ with respect to a $\check{G}$-local
system $E$. Then for the diagram 
$$
\Bun_G\times X^{(d)}\;\getsup{\gp\times\supp}\;
\cH^{+,d}_G\toup{\gq} \Bun_G
$$ 
and any local system $W$ on $X$ we have
$$
(\gp\times \supp)_!(\gq^*K\otimes\cL^d_W)[-d_G](\frac{-d_G}{2})\iso K\boxtimes (W\otimes
V^{\gamma}_E)^{(d)}[d](\frac{d}{2})
\eqno{\square}
$$
\end{Pp}

\smallskip
\noindent
2.6.1 \  For a $T$-torsor $\cF_T$ on $X$ denote by $\Bunb_N^{\cF_T}$
the stack of collections $(\cF_G,\kappa)$, where $\cF_G\in\Bun_G$ and
for $\check{\lambda}\in\check{\Lambda}^+$
$$
\kappa^{\check{\lambda}}: \cL^{\check{\lambda}}_{\cF_T}\hook{}\cV^{\check{\lambda}}_{\cF_G}
$$
are inclusions of coherent sheaves on $X$ satisfying Pl\" ucker relations (as in \cite{FGV}, 
section 2.1.2). 
 
 The open substack $j: \Bun_N^{\cF_T}\hook{}\Bunb_N^{\cF_T}$ is given by the condition that
all $\kappa^{\check{\lambda}}$ are maximal embeddings. 

 Consider the stack of
pairs $(\cF_T,\tilde\omega)$, where $\cF_T$ is a $T$-torsor on $X$ and $\tilde\omega$ is
a trivial conductor, that is, $\tilde\omega$ is a collection of isomorphisms
$\tilde\omega_i:\cL^{\check{\alpha}_i}_{\cF_T}\,\iso\,\Omega$ for each
$i\in\cal I$. The exact sequence $1\to Z(G)\to T\to \prod_{i\in\cal I}\Gm\to 1$,
where the second map is $\prod_{i\in\cal I} \check{\alpha}_i$, shows that this stack is
noncanonically isomorphic to $\Bun_{Z(G)}$ (recall that by our assumption
$Z(G)$ is connected).

Fix a section $T\to B$. Then for each pair $(\cF_T, \tilde\omega)$ we have 
the evaluation map $\ev^{\tilde\omega} :\Bun_N^{\cF_T}\to\A^1$ (cf. \cite{FGV}, section
4.1.1).  Fix a $T$-torsor on $X$ with a trivial conductor
$(\cF_T,\tilde\omega)$.

\begin{Rem} 
\label{Rem1}
If $(\cF'_T,\tilde\omega')$ is another $T$-torsor with trivial
conductor on $X$ then there exists a $Z(G)$-torsor $\cF_{Z(G)}$ on $X$ and an
isomorphism $\cF_T\otimes\cF_{Z(G)}\iso\cF'_T$ with the following property. Let
$\Bun_N^{\cF_T}\to \Bun_N^{\cF'_T}$ be the isomorphism that sends $\cF_B$ to
$(\cF_B\times\cF_{Z(G)})/Z(G)$, where $Z(G)$ acts diagonally. Then the diagram commutes
$$
\begin{array}{ccc}
\Bun_N^{\cF_T} & \iso &\Bun_N^{\cF'_T}\\
& \searrow\lefteqn{\scriptstyle \ev^{\tilde\omega}} & \downarrow\lefteqn{\scriptstyle
\ev^{\tilde\omega'}}\\
&& \A^1
\end{array}
$$
\end{Rem}

\noindent
2.6.2 \  For $d\ge 0$ let $\cY_d$ be the stack of collections $(\cF_G, D\in X^{(d)},
\kappa)$, where
$(\cF_G,\kappa)\in \Bunb_N^{\cF'_T}$ with $\cF'_T=\cF_T(w_0(\gamma)D)$. So, for each
$\check{\lambda}\in\check{\Lambda}_S^+$ we have an embedding of coherent sheaves
$$
\kappa^{\check{\lambda}}: \cL^{\check{\lambda}}_{\cF_T}\hook{}\cV^{\check{\lambda}}_{\cF_G},
$$
and $\kappa^{\check{\omega}_0}$ induces an isomorphism
$\cL^{\check{\omega}_0}_{\cF_T}(D)\,\iso\,
\cV^{\check{\omega}_0}_{\cF_G}$. 

 Let $_0\cY_d\subset \cY_d$ denote the open substack given by the condition
$(\cF_G,\kappa)\in\Bun_N^{\cF'_T}$. 

 Consider the fibred product $\cY_d\times_{\Bun_G}\cH^{+,k}_G$, where the map
$\cH^{+,k}_G\to\Bun_G$ is $\gp$. For $k\ge 0$ we have a proper representable map
$$
\gq_{\cY}:\cY_d\times_{\Bun_G}\cH^{+,k}_G\to \cY_{d+k}
$$
that sends $(\cF_G, \kappa, \cF'_G, \beta: \cF_G\mid_{X-D'}\iso
\cF'_G\mid_{X-D'})$ to $(\cF'_G, \kappa')$, where ${\kappa'}^{\check{\lambda}}$ are 
the compositions
$$
\cL^{\check{\lambda}}_{\cF_T}\hook{\kappa^{\check{\lambda}}}\cV^{\check{\lambda}}_{\cF_G}
\hook{} \cV^{\check{\lambda}}_{\cF'_G}
$$

 Given a local system $W$ on $X$, define the sheaf $\cP^d_{W,\psi}$ on $\cY_d$ as follows.
Consider the open immersion  $j:\Bun_N^{\cF_T}\to \Bunb_N^{\cF_T}\iso\cY_0$. Set
$$
\cP^0_{W,\psi}=\cP^0_{\psi}=j_!(\ev^{\tilde\omega})^*\cL_{\psi}[d_N](\frac{d_N}{2}),
$$
where $d_N=\dim\Bun_N^{\cF_T}$. 
By (Theorem 2, \cite{FGV}), this is a perverse sheaf and
$\DD(\cP^0_{\psi})\iso \cP^0_{\psi^{-1}}$. For $d>0$ set 
$$
\cP^d_{W,\psi}=\gq_{\cY!}(\cP^0_{\psi}\boxtimes \cL^d_W)[-d_G](\frac{-d_G}{2})
$$
It is easy to see that $\Qlb\boxtimes\cL^d_W$ is ULA with respect to the projection
$\cY_0\times_{\Bun_G}\cH^{d,+}_G\to\cY_0$. So, by (property 5, sect. 5.1.2 \cite{BG}), 
$$
\DD(\cP^0_{\psi}\boxtimes \cL^d_W)[-d_G](\frac{-d_G}{2}))\;\iso\;
\cP^0_{\psi^{-1}}\boxtimes \cL^d_{W^*}[-d_G](\frac{-d_G}{2})
$$
Therefore, $\DD(\cP^d_{W,\psi})\,\iso\, \cP^d_{W^*,\psi^{-1}}$ canonically.

Let $^{rss}\cY_d\subset {_0\cY_d}$ be the preimage of $^{rss}X^{(d)}$ under
the projection $_0\cY_d\to X^{(d)}$.

\begin{Pp}
\label{Pp_irreducibility}
For any local system $W$ on $X$, $\cP^d_{W,\psi}$ is a perverse sheaf on
$\cY_d$, the Goresky-MacPherson extension from $^{rss}\cY_d$. 
If $W$ is irreducible then $\cP^d_{W,\psi}$ is irreducible.
\end{Pp} 
The proof is found in Section 3.2.

\bigskip
\noindent
2.7 \  Let $\pi^0: \Bun_N^{\cF_T}\to \Bun_G$ be the projection. 
\begin{Def} Let $K$ be a $E$-Hecke
eigensheaf on $\Bun_G$, where $E$ is a $\check{G}$-local system on $X$. 
We say that $K$ is \select{generic normalized} if it is equiped with an isomorphism
$$
\RG_c(\Bun_N^{\cF_T},\; \cP^0_{\psi}\otimes
\pi^{0*}K)\;\iso\;\Qlb[d_G](\frac{d_G}{2})
$$
By Remark~\ref{Rem1}, this property does not depend
(up to a tensoring $K$ by a 1-dimensional vector space) on our choice of the pair
$(\cF_T,\tilde\omega)$.
\end{Def}

 Write $\Bun_G^d$ for the connected component of
$\Bun_G$ given by $\deg\cF_G=\deg\cF_T+d\theta$. Let $\pi:\cY_d\to\Bun^d_G$
and
$\phi:\cY_d\to X^{(d)}$ be the projections. From Proposition~\ref{Pp2} one derives

\begin{Cor} 
\label{Cor_1}
Let $K$ be a generic normalized $E$-Hecke
eigensheaf on $\Bun_G$. Let $W$ be any local system on $X$. Then for each $d\ge 0$ one has
$$
\phi_!(\pi^*K\otimes\cP^{d}_{W,\psi})\;\iso\; (W\otimes
V^{\gamma}_E)^{(d)}[d+d_G](\frac{d+d_G}{2})\eqno{\square}
$$
\end{Cor} 

 For a $\check{G}$-local system $E$ pick a $\check{G}$-local system $E^*$ such that 
$V^{\check{\lambda}}_{E^*}\;\iso\; (V^{\check{\lambda}}_E)^*$ for all
$\check{\lambda}\in\check{\Lambda}^+$. Let $K$ be a $E$-Hecke eigensheaf then
$\DD K$ is a $E^*$-Hecke eigensheaf. Assume that $\DD K$ is generic normalized
then from Corolary~\ref{Cor_1} we get an isomorphism
$$
\phi_!(\pi^* (\DD K)\otimes\cP^d_{W,\psi})\,\iso\, (W\otimes
V^{\gamma}_{E^*})^{(d)}[d+d_G](\frac{d+d_G}{2})
$$
By adjunction, it yields a nonzero map
$$
\DD K\to \pi_*\R\HOM(\cP^d_{W,\psi}, \; \phi^!(W\otimes
V^{\gamma}_{E^*})^{(d)})[d+d_G](\frac{d+d_G}{2})
$$
Dualizing, we see that this is equivalent to providing a nonzero map
\begin{equation}
\label{map_17}
\pi_!(\cP^d_{W,\psi}\otimes \phi^*(W^*\otimes V^{\gamma}_E)^{(d)})\to 
K[d_G-d](\frac{d_G-d}{2})
\end{equation}
Set $W=V^{\gamma}_E$, so we have a canonical map $\Qlb\to (W^*\otimes
V^{\gamma}_E)^{(d)}$ on $X^{(d)}$. Composing with (\ref{map_17}) we get a morphism
\begin{equation}
\label{map_18}
\pi_!\cP^d_{W,\psi}\to K[d_G-d](\frac{d_G-d}{2})
\end{equation}

\begin{Con}[geometric Langlands] 
\label{Con_Langlands} 
Let $E$ be a $\check{G}$-local system on $X$.
Assume that
$W=V^{\gamma}_E$ is irreducible and satisfies the condition

\medskip

{\rm (A)} If $E'$ is a $\check{G}$-local system on $X$ such that
$V^{\gamma}_{E'}\,\iso\, V^{\gamma}_E$ then $E'\,\iso\,
E$.

\medskip
\noindent
Then there exists $N>0$ and for each $d\ge N$ a nonempty open substack $U_d\subset \Bun^d_G$
with the following property. There exists a $E$-Hecke eigensheaf $K$ on $\Bun_G$ such
that
\begin{itemize}
\item both $K$ and $\DD K$ are generic normalized;
\item for $d\ge N$ the complex $\pi_!\cP^d_{W,\psi}\mid_{U_d}$ is placed in perverse
degrees $\le d-d_G$, and the map (\ref{map_18}) induces an isomorphism 
$$
\cH^{d-d_G}(\pi_!\cP^d_{W,\psi})\;\iso\; K\mid_{U_d}
$$ 
on the top perverse
cohomology sheaves; 
\item $K$ is an irreducible perverse sheaf over each $\Bun^d_G$, which does not vanish 
over $U_d$.
\end{itemize}
\end{Con}

\begin{Rems} i) The sheaf $K$ from Conjecture~\ref{Con_Langlands} is unique up to an
isomorphism if it exists. \\
ii) For any local system $W$ on $X$ we have
$\pi_!\cP^d_{W,\psi}\;\iso\;\Av^d_W(\pi_!\cP^0_{\psi})$ naturally.
\end{Rems}

\medskip\noindent
2.8 {\scshape Informal motivation \ } If the ground field $k$ was finite then
according to Langlands' spectral decomposition theorem (\cite{MW}), each
function from $L^2(\Bun_G(k))$ would be written as linear combination (more precisely, 
a direct integral) of Hecke eigenfunctions.

 Conjecturally, some version of spectral decomposition should exist for the derived category
$\D(\Bun_G)$ itself. We also have an analog of the scalar product of two objects $K_1,K_2\in
\D(\Bun_G)$, which is the cohomology $\RG_c(\Bun_G, K_1\otimes\DD(K_2))$ (we ignore here all
convergence questions). 

 Let $E$ be a $\check{G}$-local system on $X$ satisfying the assumptions of
Conjecture~\ref{Con_Langlands}. One may hope that to $E$ is associated a
$E$-Hecke eigensheaf $K$, which is unique in appropriate sense.

 Since $K$ is expected to be generic normalized, the "scalar product" of $\pi_!\cP^0_{\psi}$
and $K$ should equal "one". That is, $K$ should appear in the spectral decomposition of 
$\pi_!\cP^0_{\psi}$ with multiplicity one.  By Proposition~\ref{Pp2},
the functor $\Av^d_W$ applied to $\pi_!\cP^0_{\psi}$ with $d$ large enough, will kill all the
terms in the spectral decomposition of $\pi_!\cP^0_{\psi}$ except $K$ itself. So, roughly
speaking, $\Av^d_W(\pi_!\cP^0_{\psi})$ should equal $K$ tensored by some constant complex. 

\bigskip\noindent
2.9 \ {\scshape Stratifications} \  For $\mu\in \Lambda^{pos}$ denote by $X^{\mu}$ the moduli
scheme of $\Lambda^{pos}$-valued divisors of degree $\mu$. If $\mu=\sum_{i\in\cal I}
a_i\alpha_i$ then $X^{\mu}=\prod_{i\in\cal I} X^{(a_i)}$.

\smallskip

 For $D\in X^{(d)}$ consider $\cF'_T=\cF_T(w_0(\gamma)D)$. The
stack $\cY_d$ is the stack of pairs: $D\in X^{(d)}$ and a point $(\cF_G,\kappa)\in
\Bunb_N^{\cF'_T}$. Recall that $\Bunb_N^{\cF'_T}$ is stratified by locally closed substacks
$_{\mu}\Bun_N^{\cF'_T}$ indexed by $\mu\in\Lambda^{pos}$. Namely, $(\cF_G,\kappa)$ lies in 
$_{\mu}\Bun_N^{\cF'_T}$ iff there exists a divisor $D^{pos}\in X^{\mu}$ such that for all
$\check{\lambda}\in\check{\Lambda}^+$ the meromorphic maps
$$
\cL^{\check{\lambda}}_{\cF'_T}(\<D^{pos}, \check{\lambda}\>)\to \cV^{\check{\lambda}}_{\cF_G}
$$
are regular everywhere and maximal.
We have a projection
$$
_\mu\Bun_N^{\cF'_T}\to X^{\mu}
$$ 
whose fibre over $D^{pos}$ is isomorphic to
$\Bun_N^{\tilde\cF_T}$, where $\tilde\cF_T\;\iso\; \cF_T(w_0(\gamma)D+D^{pos})$.

 Denote by $_\mu\cY_d\subset \cY_d$ the locally closed substack given by 
$(\cF_G,\kappa)\in {_\mu\Bun_N^{\cF'_T}}$. Let
$$
\phi^{\mu}: {_\mu\cY_d}\to X^{(d)}\times X^{\mu}
$$
be the projection. Let $X^{d,\mu}\subset X^{(d)}\times X^{\mu}$ be the closed
subscheme given by the condition: $\gamma D+w_0(D^{pos})$ is dominant. This condition
ensures that the maps $\tilde\omega_i:\cL^{\check{\alpha}_i}_{\tilde\cF_T}\to \Omega$ are
regular for all $i\in\cal I$. Let $_\mu\cY_d^+\subset {_\mu\cY_d}$ be the preimage of
$X^{d,\mu}$ under $\phi^{\mu}$. So, we have the evaluation map 
$$
\ev_{\mu}:{_\mu\cY_d^+}\to
\A^1
$$

 Note that if $(D,D^{pos})\in X^{d,\mu}$ then $\gamma D+w_0(D^{pos})$ is a
$\Lambda^+_{G,S}$-valued divisor of degree $d\gamma+w_0(\mu)$. In this way
$X^{d,\mu}$ is the moduli scheme of $\Lambda^+_{G,S}$-valued divisors on $X$ of degree
$d\gamma+w_0(\mu)$.

\bigskip\noindent
2.10 \  For any local system $W$ on $X$ let $\tilde\cP^d_{W,\psi}$ be the complex obtained by
replacing in the definition of $\cP^d_{W,\psi}$ Laumon's sheaf by Springer's sheaf 
$$
\tilde\cP^d_{W,\psi}=\gq_{\cY !}(\cP^0_{\psi}\boxtimes\Spr^d_W)[-d_G](\frac{-d_G}{2})
$$
\begin{Pp} 
\label{Pp_3}
1) For each $\mu\in\Lambda^{pos}$ the restriction of $\tilde\cP^d_{W,\psi}$ to
$_\mu\cY_d$ is supported by $_\mu\cY^+_d$ and is isomorphic to 
$$
\phi^{\mu *}\tilde W^{d,\mu}\otimes
\ev_{\mu}^*\cL_{\psi}\otimes(\Qlb(\frac{1}{2})[1])^{\otimes d+d_N+\<\gamma d-\mu, \,
2\check{\rho}\>}
$$
for some sheaf $\tilde W^{d,\mu}$ on $X^{d,\mu}$. Here $\tilde W^{d,\mu}$ is placed in usual
degree zero. Let $W^{d,\mu}$ denote the corresponding sheaves for $\cP^d_{W,\psi}$.\\
2) If $(D=\sum d_kx_k, \; D^{pos}=\sum \mu_kx_k)$ is a $k$-point of $X^{d,\mu}$ then
the fibre of $W^{d,\mu}$ at this point is the tensor product over all $x_k$ 
$$
\mathop{\otimes}\limits_{k}\; (\mathop{\oplus}\limits_{\nu} \Hom(V^{d_k\gamma+w_0(\mu_k)},
(V^{\gamma})_{\nu})\otimes (W_{x_k})_{\nu}),
$$
the inside sum being taken over partitions $\nu$ of $d_k$ of length $\le r$.
\end{Pp}

 The proof is found in Section 3.1.
Proposition~\ref{Pp_3} together with 
Corolary~\ref{Cor_1} suggest the following 
conjecture.
\begin{Con}
Let $D=\sum \lambda_kx_k$ be a divisor on $X$ with $\lambda_k$ dominant coweights, let
$\lambda$ be the degree of $D$. Denote by
$a:\Bun_N^{\tilde\cF_T}\to\Bun_G$ the projection, where
$\tilde\cF_T=\cF_T(w_0(D))$. Let $\ev_{\lambda}: \Bun_N^{\tilde\cF_T}\to\A^1$ be the evaluation
map given by the conductor data. 
Let $E$ be a $\check{G}$-local system
on $X$, $K$ be a generic normalzed $E$-Hecke eigensheaf on $\Bun_G$. Then
$$
\RG_c(\Bun_N^{\tilde\cF_T}, a^*K\otimes \ev_{\lambda}^*\cL_{\psi})\otimes
(\Qlb(\frac{1}{2})[1])^{d_N-d_G+\<\lambda, 2\check{\rho}\>}\;\,\iso\;\,
\mathop{\otimes}\limits_{k} (V^{\lambda_k}_E)_{x_k}
$$
\end{Con}

\noindent
2.11 \ Consider the diagram
$$
\cY_d\times X\;\getsup{\gp_{\cY}\times\supp}\;
\cY_d\times_{\Bun_G}\cH^{+,1}_G \toup{\gq_{\cY}}\cY_{d+1}
$$
where $\gp_{\cY}: \cY_d\times_{\Bun_G}\cH^{+,1}_G\to \cY_d$ is the projection.
For any local system $W$ on $X$ define a natural map
\begin{equation}
\label{map_1}
(\gp_{\cY}\times\supp)_!\gq_{\cY}^*\cP^{d+1}_{W,\psi}
\otimes\Qlb(\frac{1}{2})[1]^{\otimes \<\gamma, 2\check{\rho}\>}
\to \cP^d_{W,\psi}\boxtimes W \otimes\Qlb(\frac{1}{2})[1]
\end{equation}
as follows.
Consider the diagram
$$
\begin{array}{ccc}
\tilde\cH^{+,d}_G\times_{\Bun_G}\cH^{+,1}_G  & \iso & \tilde\cH^{+,\, d+1}_G\\
\downarrow\lefteqn{\scriptstyle p\times\id} && \downarrow\lefteqn{\scriptstyle p}\\
\cH^{+,d}_G\times_{\Bun_G}\cH^{+,1}_G & \toup{\conv} & \cH^{+,\, d+1}_G
\end{array}
$$
Clearly, $\cL^{d+1}_W$ is a direct summand of $\conv_!(\cL^d_W\boxtimes\supp^*W)\otimes
\Qlb(\frac{1}{2})[1]^{\otimes 1+\<\gamma, 2\check{\rho}\>}$. From the diagram
$$
\begin{array}{ccc}
\cY_0\times_{\Bun_G}\cH^{+,d}_G\times_{\Bun_G}\cH^{+,1}_G & \toup{\id\times\conv} &
\cY_0\times_{\Bun_G}\cH^{+,d+1}_G\\
\downarrow && \downarrow\\
\cY_d\times_{\Bun_G}\cH^{+,1}_G & \toup{\gq_{\cY}} & \cY_{d+1}
\end{array}
$$
we see that $\cP^{d+1}_{W,\psi}$ is a direct summand of 
$$
\gq_{\cY *}(\cP^d_{W,\psi}\boxtimes
\supp^* W)\otimes\Qlb(\frac{1}{2})[1]^{\otimes 1+\<\gamma, 2\check{\rho}\>}
$$
This yields a morphism
$$
\gq_{\cY}^*\cP^{d+1}_{W,\psi}\to (\cP^d_{W,\psi}\boxtimes
\supp^* W)\otimes\Qlb(\frac{1}{2})[1]^{\otimes 1+\<\gamma, 2\check{\rho}\>}
$$
Since $\gp_{\cY}\times\supp: \cY_d\times_{\Bun_G}\cH^{+,1}_G\to \cY_d\times X$ is smooth
of relative dimension $\<\gamma, 2\check{\rho}\>$, we get a map
$$
\gq_{\cY}^*\cP^{d+1}_{W,\psi}\to (\gp_{\cY}\times\supp)^!(\cP^d_{W,\psi}\boxtimes W)
\otimes\Qlb(\frac{1}{2})[1]^{\otimes 1-\<\gamma, 2\check{\rho}\>}
$$
and, by adjunction, the desired map (\ref{map_1}).

\begin{Pp} 
\label{Pp_on_Hecke}
1) The complex 
\begin{equation}
\label{complex_1}
(\gp_{\cY}\times\supp)_!\gq_{\cY}^*\cP^{d+1}_{W,\psi}\otimes\Qlb(\frac{1}{2})[1]^{\otimes
\<\gamma, 2\check{\rho}\>}
\end{equation}
is placed in perverse degrees $\le 0$.\\
2) For any $\mu\in \Lambda^{pos}$, the restriction of {\rm (\ref{complex_1})} to 
$_{\mu}\cY_d\times X$ is supported by
$_{\mu}\cY_d^+\times X$ and is isomorphic to the tensor product of
$$
\ev_{\mu}^*\cL_{\psi}\otimes (\Qlb(\frac{1}{2})[1])^{\otimes d+1+d_N+\<\gamma d-\mu, \,
2\check{\rho}\>}
$$
with some sheaves $W_0^{d,\mu}$ coming from $X^{d,\mu}\times X$. Here $W_0^{d,\mu}$ is
placed in usual degree zero.
\end{Pp}

 The proof is given in Section 3.3

\begin{Rems}
i) One may show that for any $\mu\in \Lambda^{pos}$ the restriction of (\ref{map_1}) to
$_\mu\cY_d\times X$ comes from a morphism of sheaves $W_0^{d,\mu}\to W^{d,\mu}\boxtimes W$ on
$X^{d,\mu}\times X$. The latter map is an isomorphism over the open substack of
$X^{d,\mu}\times X$ classifying triples
$(D,D^{pos},x)$ such that $x$ does not appear in $D$. Therefore, (\ref{map_1}) is an
isomorphism over the locus of $(\cF,\kappa,D,x)\in\cY_d\times X$ such that $x$ does not appear
in $D$.\\  
ii) In the situation of Conjecture~\ref{Con_Langlands} we expect that the map
(\ref{map_1}) yields the Hecke property of $K$ corresponding to the coweight $\gamma$.
Moreover, it  should also yield the Hecke properties corresponding to all
$\lambda\in\Lambda^+_{G,S}$ (as it indeed happens for $\GL_n$). Define the Hecke functor
$\H_{\cY}:\D(\cY_{d+1})\to
\D(\cY_d\times X)$ by 
$$
\H_{\cY}(F)=(\gp_{\cY}\times\supp)_!\gq_{\cY}^*F\otimes\Qlb(\frac{1}{2})[1]^{\otimes
\<\gamma, 2\check{\rho}\>}
$$
Note that the cohomological shifts in the definition of
the Hecke functor $\H:\D(\Bun^{d+1}_G)\to\D(\Bun^d_G\times X)$ corresponding to $\gamma$
and $\H_{\cY}$ differ by one! So, the Hecke property of $K$ can not be simply the
push-forward of (\ref{map_1}) with respect to $\pi: \cY_d\to\Bun_G$.
\end{Rems}

\medskip\noindent
2.12 \  Let $\omega$ be a generator of the group of coweights orthogonal to all roots. Since
the image of $\omega$ in $\pi_1(G)$ is not zero, we assume that this image
equals $d_{\omega}\theta$ for some $d_{\omega}>0$. Since $d_{\omega}\gamma-\omega$ is dominant,
we have $d_{\omega}\gamma-\omega\in\Lambda^{pos}$ and $\omega\in \Lambda^+_{G,S}$.
Consider the map $\gq^{\omega}:\Bun_G\times X\to \Bun_G$ sending $(\cF_G,x)$ to $\cF'_G$, 
where $\cF'_G$ and $\cF_G$ are identified over $X-x$ and
$$
\cV^{\check{\lambda}}_{\cF_G}(\<\omega,\check{\lambda}\>x)\;\iso\;
\cV^{\check{\lambda}}_{\cF'_G}
$$
for each $\check{\lambda}\in\check{\Lambda}^+$. Let also 
$$
\gq_{\cY}^{\omega}: \cY_d\times X\to
\cY_{d+d_{\omega}}
$$ 
be the map sending $(\cF_G,\kappa,x)$ to $(\cF'_G,\kappa')$, where
$\cF'_G=\gq^{\omega}(\cF_G,x)$ and $\kappa^{\prime\check{\lambda}}$ is the composition
$$
\cL^{\check{\lambda}}_{\cF_T}\hook{\kappa^{\check{\lambda}}} \cV^{\check{\lambda}}_{\cF_G}
\hook{} \cV^{\check{\lambda}}_{\cF'_G}
$$
for all $\check{\lambda}\in\check{\Lambda}^+_S$. Let $E$ be a $\check{G}$-local
system on $X$ and set $W=V^{\gamma}_E$. Then there is a natural map
\begin{equation}
\label{map_2}
(\gq^{\omega}_{\cY})^*\cP^{d+d_{\omega}}_{W,\psi}\otimes\Qlb(\frac{1}{2})[1]^{\otimes
1-d_{\omega}}
\to  \cP^d_{W,\psi}\boxtimes
V^{\omega}_E\otimes\Qlb(\frac{1}{2})[1]
\end{equation}
This is not an isomorphism in general, and one may show that the LHS of (\ref{map_2}) is placed
in perverse degrees $\le 0$.

\medskip\noindent
2.13 {\scshape Whittaker sheaves\ } Let $\K(\cY_d)$ denote the Grothendieck ring of the
triangulated category $\D(\cY_d)$. 

To each $\check{G}$-local system $E$ on $X$ and $d\ge 0$ we attach \select{the Whittaker
sheaf} $\cW^d_{E,\psi}\in \K(\cY_d)$ defined as follows. 

 Let $\mu\in\Lambda^{pos}$ be such that $d\gamma+w_0(\mu)$ is
dominant. Let $\tau$ be a partition of $d\gamma+w_0(\mu)$, that is, a way to
write $d\gamma+w_0(\mu)=\sum_i n_i\lambda_i$ with $\lambda_i\in\Lambda^+_{G,S}$ pairwise
different and $n_i>0$. Let $^{\tau}X\subset\prod_i X^{(n_i)}$ be the complement to all
diagonals. We consider $^{\tau}X\subset X^{d,\mu}$ as the locally-closed subscheme
classifying divisors $\sum_k \lambda_k x_k$ of degree $d\gamma+w_0(\mu)$
with $x_k\in X$ pairwise different. 

 Let $^{\tau}E$ denote the restriction of 
$$
\mathop{\boxtimes}\limits_i
(V^{\lambda_i}_E)^{(n_i)}
$$ 
under $^{\tau}X\hook{}\prod_i X^{(n_i)}$. Let $^{\tau}\cY\subset {_{\mu}\cY_d^+}$
be the preimage of $^{\tau}X$ under $\phi^{\mu}: {_{\mu}\cY_d^+}\to X^{d,\mu}$. 

\begin{Def} Set $\cW^d_{E,\psi}\in \K(\cY_d)$ to be the (unique) complex with the
following properties. Its $*$-restriction to each stratum $_{\mu}\cY_d$
is supported by $_{\mu}\cY_d^+$. If $d\gamma+w_0(\mu)$ is dominant then
for any partition $\tau$ of $d\gamma+w_0(\mu)$ the
$*$-restriction of $\cW^d_{E,\psi}$ to $^{\tau}\cY$ is
$$
\phi^{\mu *} (^{\tau}E)\otimes\ev_{\mu}^*\cL_{\psi}\otimes \Qlb(\frac{1}{2})[1]^{\otimes
d+d_N+\<d\gamma-\mu, 2\check{\rho}\>}
$$
\end{Def} 

 The sheaf $\cW^d_{E,\psi}$ should satisfy the Hecke property, in particular
we suggest

\begin{Con} Recall the diagram (cf. Sect. 2.11)
$$
\cY_d\times X\;\getsup{\gp_{\cY}\times\supp}\;
\cY_d\times_{\Bun_G}\cH^{+,1}_G \toup{\gq_{\cY}}\cY_{d+1}
$$
There is a canonical isomorphism in the Grothendieck ring $K(\cY_d\times X)$
$$
(\gp_{\cY}\times\supp)_!\gq_{\cY}^*\cW^{d+1}_{E,\psi}
\otimes\Qlb(\frac{1}{2})[1]^{\otimes \<\gamma, 2\check{\rho}\>}
\to \cW^d_{E,\psi}\boxtimes V^{\gamma}_E \otimes\Qlb(\frac{1}{2})[1]
$$
\end{Con}

\medskip\noindent
2.13.1 We don't know if $\Lambda^+_{G,S}$ is a free semigroup in general, however this
is the case for our examples $\GL_n, \GSp_{2n}, \GSpin_{2n+1}$ (cf. appendix). 

\smallskip

 Assuming that $\Lambda^+_{G,S}$ is a free semigroup, we can describe $\cW^d_{E,\psi}$
more precisely, namely "glue" the pieces on the strata $^{\tau}\cY$ to get a sheaf on
$_{\mu}\cY_d^+$. To do so, we will glue the sheaves $^{\tau}E$ to get a constructible
sheaf $AE^{d,\mu}$ on $X^{d,\mu}$ (here `A' stands for `automorphic').

\smallskip

 Let $\lambda_1,\ldots,\lambda_m$ be free generators of $\Lambda^+_{G,S}$ thus
yielding $\Lambda^+_{G,S}\,\iso\, (\ZZ_+)^m$. Given $d\ge 0$ and $\mu\in\Lambda^{pos}$
with $d\gamma+w_0(\mu)=\sum_{i=1}^m a_i\lambda_i$ dominant, we get 
$$
X^{d,\mu}\,\iso\, \prod_{i=1}^m
X^{(a_i)}
$$
Consider the sheaf 
\begin{equation}
\label{sheaf_divisors_Whittaker}
\mathop{\boxtimes}\limits_{i=1}^m (V^{\lambda_i}_E)^{(a_i)}
\end{equation}
on $X^{d,\mu}$. Let $D=\sum_k \nu_k x_k$ be a $k$-point of $X^{d,\mu}$, where $x_k$ are
pairwise different and $\nu_k\in\Lambda^+_{G,S}$. Write $\nu_k=\sum_{i=1}^m
a_{i,k}\lambda_i$ for each $k$. The fibre of (\ref{sheaf_divisors_Whittaker}) at $D$ is
$$
\otimes_{x_k} \otimes_{i=1}^m
\Sym^{a_{i,k}}(V^{\lambda_i}_E)_{x_k}
$$
There is (a unique up to a nonzero multiple) inclusion of $\check{G}$-modules
$V^{\nu_k}\hook{} \otimes_{i=1}^m \Sym^{a_{i,k}}(V^{\lambda_i})$. This yields a map
\begin{equation}
\label{fibre_of_AE}
\otimes_{x_k} (V^{\nu_k}_E)_{x_k}\hook{} \otimes_{x_k} \otimes_{i=1}^m
\Sym^{a_{i,k}}(V^{\lambda_i}_E)_{x_k}
\end{equation}

 The following is borrowed from \cite{Ly2}.

\begin{Pp} Assume that $\Lambda^+_{G,S}$ is a free semigroup.
There is a unique constructible subsheaf $AE^{d,\mu}\subset 
\mathop{\boxtimes}\limits_{i=1}^m (V^{\lambda_i}_E)^{(a_i)}$ whose fibre at 
any $k$-point $D=\sum_k \nu_k x_k$ of $X^{d,\mu}$ is the image of (\ref{fibre_of_AE}).
\QED
\end{Pp}

\bigskip\bigskip
\centerline{{\scshape 3. Some proofs}}

\bigskip\noindent
3.1 \  Recall that $\Gr_G$ is stratified by locally closed ind-subschemes $S^{\mu}$ indexed by
all coweights $\mu\in\Lambda$. Informally, $S^{\mu}$ is the $N(\hat K)$-orbit of the point
$\mu(t)\in\Gr_G$, where $\hat K=k((t))$. We refer the reader to \cite{FGV}, Section~7.1 for the
precise definition. 

 Recall the following notion from \select{loc.cit.}, section 7.1.4. Set $\hat\cO=k[[t]]$. Let
$\hat\Omega$ denote the completed module of relative differentials of
$\hat\cO$ over $k$ (so, $\hat\Omega$ is a free $\hat\cO$-module generated by $dt$).
Given a coweight $\eta\in\Lambda$ and isomorphisms
$$
s_i:\hat\cO(\<\eta,\check{\alpha}_i\>)\,\iso\,\hat\Omega
$$ 
for each $i\in\I$, one defines an admissible
character $\chi_{\eta}: N(\hat K)\to\Ga$ of conductor $\eta$
as the sum $\chi_{\eta}=\sum_{i\in\I} \chi^i$,
where $\chi^i: N(\hat K)\to\Ga$ is the composition
$$
N(\hat K)\to N/[N,N](\hat K)\toup{u_i}\Ga(\hat K)\toup{s_i}\hat\Omega(\hat K)\toup{\Res}\Ga
$$
Here $u_i:N/[N,N]\to\Ga$ is the natural coordinate corresponding to the simple root
$\check{\alpha}_i$. By (\select{loc.cit.}, Lemma~7.1.5), for $\nu\in\Lambda$ there exists a 
$(N(\hat K),\chi_{\eta})$-equivariant function $\chi^{\nu}_{\eta}: S^{\nu}\to\A^1$ if and only
if $\nu+\eta\in\Lambda^+$. In the latter case this function is unique up to an additive
constant. 

\bigskip

\begin{Prf}\select{of Proposition~\ref{Pp_3}\ } 
Let $(\sum d_kx_k,\, \cF_G)$ be a $k$-point of $X^{(d)}\times\Bun_G$.
The fibre of $\supp\times\gq:\cH^{+,d}_G\to X^{(d)}\times\Bun_G$ over
this point identifies with 
$$
\prod_k \Grb^{d_k\gamma}_G
$$
The restriction of $\Spr^d_W$ to this fibre is the tensor product of $(\otimes_k
(W_{x_k}^{\otimes d_k}))[d+d_G](\frac{d+d_G}{2})$ with
$$
\mathop{\boxtimes}\limits_k \;\,(
\mathop{\oplus}\limits_{\nu\le d_k\gamma} \cA_{\nu}\otimes V_{\nu,k})
$$
where the inside sum is taken over dominant coweights $\nu\in\Lambda^+$ such that $\nu\le
d_k\gamma$, and $V_{\nu,k}$ are some vector spaces. Let 
$$
^0\gq_{\cY}:\Bun_N^{\cF_T}\times_{\Bun_G}\cH^{+,d}_G\to \cY_d
$$ 
denote the restriction of
$\gq_{\cY}$ to the open substack $\Bun_N^{\cF_T}\times_{\Bun_G}\cH^{+,d}_G\subset
\cY_0\times_{\Bun_G}\cH^{+,d}_G$. 

 Fix a $k$-point $y\in\cY_d$, it is given by $(D=\sum d_kx_k,\; D^{pos}=\sum \mu_k
x_k,\; \cF_G,\kappa)$ such that $x_k$ are pairwise different (and some of $d_k$
may be zero). Let $K_y$ denote the fibre at $y$ of 
$$
^0\gq_{\cY!}(\ev^*\cL_{\psi}\boxtimes\Spr^d_W)[d_N-d_G](\frac{d_N-d_G}{2})
$$
The fibre of $^0\gq_{\cY}$ over $y$ identifies with
$$
\prod_k S^{w_0(\gamma)d_k+\mu_k}\cap \Grb^{d_k\gamma}_G
$$
For each $k$ set $\nu_k=-w_0(\gamma)d_k-\mu_k$. An equivariance argument (as in \cite{FGV},
Lemma~6.2.8) shows that $K_y$ vanishes unless all $\nu_k$ are dominant.

 By (Lemma~7.2.7(2),
\cite{FGV}) the restriction of the map
$$
\Bun_N^{\cF_T}\times_{\Bun_G}\cH^{+,d}_G\to \Bun_N^{\cF_T}\toup{\ev}\A^1
$$
to $(^0\gq_{\cY})^{-1}(y)$ becomes the sum over all $k$ of
$$
\chi^{-\nu_k}_{\nu_k}:S^{-\nu_k}\cap \Grb^{d_k\gamma} \to\A^1
$$
plus $\ev_{\mu}(\cF_G,\kappa, D, D^{pos})\in \A^1$.
By (Theorem~1, \cite{FGV}), for any $\nu\in\Lambda^+$ such that $\nu\le d_k\gamma$ the
complex
$$
\RG_c(S^{-\nu_k}\cap \Grb^{d_k\gamma}, \cA_{\nu}\otimes(\chi^{-\nu_k}_{\nu_k})^*\cL_{\psi})
$$
vanishes unless $\nu=-w_0(\nu_k)$. In the latter case, it is canonically
$\Qlb[\<\nu_k,2\check{\rho}\>](\<\nu_k,\check{\rho}\>)$. 

\smallskip

 The above equivariance argument shows also that the restriction of
$\tilde\cP^d_{W,\psi}$ to 
$_\mu\cY_d^+$, after tensoring by $\ev_{\mu}^*\cL_{\psi^{-1}}$, descends with respect to the
projection $_\mu\cY_d^+\to X^{d,\mu}$. Combining this with
Proposition~\ref{Pp_description_Laumon}, one finishes the proof of Proposition~\ref{Pp_3}. 
\end{Prf}

\medskip

 The above proof combined with (Proposition~3.2.6, \cite{BG}) also gives the following

\begin{Cor} Over the open substack $_0\cY_d\subset \cY_d$, the map $^0\gq_{\cY}:
\Bun_N^{\cF_T}\times_{\Bun_G}\cH^{+,d}_G\to \cY_d$ is an isomorphism.
\end{Cor}

\medskip\noindent
3.2 \  In this subsection we prove Proposition~\ref{Pp_irreducibility}. 

  Given a pair $d\ge 0$, $\mu\in
\Lambda^{pos}$ such that $d\gamma+w_0(\mu)$ is dominant, a \select{partition} $\tau$ of
$d\gamma+\omega_0(\mu)$ is a presentation of
$d\gamma+\omega_0(\mu)$ as a sum of nonzero elements from $\Lambda_{G,S}^+$, that is, 
$\sum_k (d_k\gamma+\omega_0(\mu_k))=d\gamma+\omega_0(\mu)$ with $d_k\ge 0,
\mu_k\in\Lambda^{pos}$ and $d_k\gamma+\omega_0(\mu_k)$ dominant for all $k$.

 Given a partition $\tau$ of $d\gamma+\omega_0(\mu)$, consider the locally closed subscheme
$^{\tau}X\subset X^{d,\mu}$, which is the moduli scheme of divisors $\sum_k
(d_k\gamma+\omega_0(\mu_k))x_k$ with $x_k$ pairwise different. Given $\tau$, if $k$ runs 
through the set consisting of $m$ elements then $\dim {^{\tau}X}=m$. Clearly, the schemes
$^{\tau}X$ form a stratification of $X^{d,\mu}$. 

 Let $^{\tau}\cY\subset {_\mu\cY^+_d}$ be the preimage of $^{\tau}X$ in $_\mu\cY^+_d$.
Suppose that $\tau$ is of length $m$, that is, $\dim {^{\tau}X}=m$. Then 
from (Lemma~7.2.4, \cite{FGV}) it follows that
$$
\dim {^{\tau}\cY}=m+d_N+\<\gamma d-\mu, 2\check{\rho}\>
$$
In particular, we have $\dim\cY_d=d+d_N+\<d\gamma, 2\check{\rho}\>$. 
Now from Proposition~\ref{Pp_3}, we learn that the restriction of $\tilde\cP^d_{W,\psi}$ to
$^{\tau}\cY$ is placed in perverse degrees $\le 0$. Moreover, the inequality is strict unless
$\mu=0$ and $m=d$. Since $\tilde\cP^d_{W,\psi}$ is self-dual (up to replacing $W$ by $W^*$ and
$\psi$ by $\psi^{-1}$), our assertion follows. \QED

\begin{Rems} 
i) As a corolary, note that the restriction of $\cP^d_{W,\psi}$ to $_0\cY_d$ identifies
canonically with
$$
\phi^{0 *}W^{(d)}\otimes \ev_0^*\cL_{\psi}\otimes (\Qlb(\frac{1}{2})[1])^{\otimes
d+d_N+\<\gamma d, \, 2\check{\rho}\>}
$$
ii) For each pair $d\ge 0$, $\mu\in
\Lambda^{pos}$ such that $d\gamma+w_0(\mu)$ is dominant, the projection $X^{d,\mu}\to
X^{(d)}$ is a finite morphism. Let $W$ and $W'$ be any local systems on $X$. Then
the complex
\begin{equation}
\label{complex_17}
\phi_!(\cP^d_{W,\psi}\otimes \cP^d_{W',\psi^{-1}})
\end{equation}
is placed in usual cohomological degree
$-2d$. This is seen by calculating this direct image with respect to the stratification of
$\cY_d$ by $_\mu\cY_d$. For $G=\GL_n$ the complex (\ref{complex_17}) is a Rankin-Selberg
integral considered in \cite{Ly}.
\end{Rems}

\medskip\noindent
3.3 \ In this subsection we prove Proposition~\ref{Pp_on_Hecke}. 

\begin{Lm} 
\label{Lm_on_characters}
Let $\lambda\in\Lambda^+, \mu\in\Lambda$. Let $w\in W$ be any element of the Weil
group such that $\mu+w\lambda\in\Lambda^+$. The function $\chi^{w\lambda}_{\mu}:
S^{w\lambda}\to\A^1$ is constant on $S^{w\lambda}\cap\Grb^{\lambda}$ if and only if
$\mu\in\Lambda^+$.
\end{Lm}
\begin{Prf} This follows from the description of $S^{w\lambda}\cap\Grb^{\lambda}$ given
in (\cite{NP}, Lemma~5.2).
\end{Prf}

\bigskip

\begin{Prf}\select{of Proposition~\ref{Pp_on_Hecke}\ }
The fibre of $\gp\times\supp: \cH^{+,1}_G\to \Bun_G\times X$ over a $k$-point $(\cF_G,x)$
is identified with
$\Gr_G^{-w_0(\gamma)}$. This fibre is stratified by the subschemes 
\begin{equation}
\label{intersect_1}
S^{-w_0w(\gamma)}\cap
\Gr_G^{-w_0(\gamma)}
\end{equation}
indexed by $w(\gamma)\in W\gamma$. By (\cite{NP}, Lemma~5.2), (\ref{intersect_1}) is an affine
space of dimension $\<\gamma+w(\gamma), \check{\rho}\>$. 

 Consider a $k$-point of $_\mu\cY_d$ given by $(\cF_G,\kappa, D, D^{pos})$. Let
$(\cF_G,\cF'_G,\beta, x\in X)\in \cH^{+,1}_G$ correspond to a point of $S^{-w_0w(\gamma)}\cap
\Gr_G^{-w_0(\gamma)}$ for some $w(\gamma)\in W\gamma$. The image of this collection 
under the map 
$$
\gq_{\cY}:\cY_d\times_{\Bun_G}\cH^{+,1}_G\to\cY_{d+1}
$$ 
is the point of $_{\mu'}\cY_{d+1}$ given by
$$
(\cF'_G,\kappa', D+x, D^{pos}+w_0w(\gamma)x-w_0(\gamma)x)
$$
with $\mu'=\mu+w_0w(\gamma)-w_0(\gamma)$. 
If $\gamma D+w_0(D^{pos})+w(\gamma)x$ is not
dominant then, by Proposition~\ref{Pp_3}, the stratum (\ref{intersect_1}) does not
contribute to the direct image (\ref{complex_1}). 
 
 Assume $\gamma D+w_0(D^{pos})+w(\gamma)x$ dominant. 
Let $d_x\in\ZZ_+$ (resp, $\mu_x\in\Lambda^{pos}$) denote
the miltimplicity of $x$ in $D$ (resp., in $D^{pos}$).
Then the restriction of 
$\ev_{\mu'}\comp\gq_{\cY}$ to the stratum (\ref{intersect_1}) is $(N(\hat K),
\chi_{\nu})$-equivariant for some admissible character $\chi_{\nu}:N(\hat K)\to\Ga$ 
of conductor $\nu=-w_0(\gamma)d_x-\mu_x$. If $\nu$ is dominant then this restriction is 
constant and equals
$\ev_{\mu}(\cF_G,\kappa,D,D^{pos})$.
If $\nu$ is not dominant then, by
Lemma~\ref{Lm_on_characters}, the stratum (\ref{intersect_1}) does not contribute
to the direct image (\ref{complex_1}), because $\RG_c(\A^1,\cL_{\psi})=0$.

We conclude that the restriction of (\ref{complex_1}) to
$_\mu\cY_d\times X$  vanishes outside the closed substack $_\mu\cY_d^+\times X$, and is
isomorphic to 
$$
F[d+1+d_N+\<\gamma d-\mu, \,
2\check{\rho}\>]
$$ 
for some sheaf $F$ on $_\mu\cY_d^+\times X$ placed in usual degree zero. An equivariance
argument (as in the proof of Proposition~\ref{Pp_3}) assures that
$F\otimes\ev_{\mu}^*\cL_{\psi^{-1}}$ descends with respect to the projection
$_\mu\cY_d^+\times X\to X^{d,\mu}\times X$. 

 Since $\dim({_\mu\cY_d^+})\le d+d_N+\<\gamma d-\mu, 2\check{\rho}\>$ for any
$\mu\in\Lambda^{pos}$, the complex (\ref{complex_1}) is placed in perverse degrees $\le 0$.
Proposition~\ref{Pp_on_Hecke} is proved. \end{Prf}

\bigskip\bigskip
\centerline{\scshape 4. Additional structure on Levi subgroups}

\medskip\noindent
4.1 Throughout this section, we fix a standard parabolic subgroup $P\subset G$
corresponding to a subset ${\cal I}_M\subset \cal I$. Let $\Lambda_{G,P}$ be 
the quotient of $\Lambda$ by the $\ZZ$-span of $\alpha_i,
i\in{\cal I}_M$. Let 
$$
\check{\Lambda}_{G,P}=\{\check{\lambda}\in\check{\Lambda}\mid
\<\alpha_i,\check{\lambda}\>=0 \;\mbox{for}\; i\in{\cal I}_M\}
$$ 
be the dual lattice. Let $U(P)$ be the unipotent radical of $P$ and $M=P/U(P)$. Write
$\Lambda^+_M\subset\Lambda$ for the semigroup of dominant coweights for $M$, and $w_0^M$
for the longuest element of the Weil group $W_M$ of $M$. The irreducible
$\check{M}$-module of h.w. $\lambda$ is denoted $U^{\lambda}$.
Fix a section $M\to P$ of the
projection $P\to M$.

\begin{Def} Set $\Lambda^+_{M,S}=\{\lambda\in\Lambda^+_M\mid \mbox{there exists}\; w\in
W\;\mbox{such that}\; w\lambda\in \Lambda^+_{G,S}\}$.  Let $\pi_1^+(M)\subset \pi_1(M)=\Lambda_{G,P}$ denote the image of the
projection $\Lambda^+_{M,S}\to \Lambda_{G,P}$. Set 
$$
\check{\Lambda}^+_{M,S}=\{\check{\lambda}\in\check{\Lambda}^+_M\mid
\;\mbox{there exists}\; w\in W\;\,\mbox{such that}\;\,
w\check{\lambda}\in\check{\Lambda}^+_S\}
$$
\end{Def}
The following identities are straightforward.

\begin{Lm} We have
\begin{itemize}
\item[] $\Lambda^+_{M,S}=\{\lambda\in\Lambda^+_M\mid \<w(\lambda),
\check{\lambda}\>\ge 0\;\,\mbox{for all}\;\,  w\in W,
\check{\lambda}\in\check{\Lambda}^+_S\}$

\item[]
$\Lambda^+_{M,S}=\{\lambda\in\Lambda^+_M\mid \<w_0^M(\lambda),\, 
\check{\lambda}\>\ge 0\;\,
\mbox{for all}\;\, \check{\lambda}\in\check{\Lambda}^+_{M,S}\}$

\item[]
$\check{\Lambda}^+_{M,S}=\{\check{\lambda}\in\check{\Lambda}^+_M\mid
\<w(\lambda), \check{\lambda}\>\ge 0 \;\, \mbox{for all}\;\, w\in W,
\lambda\in\Lambda^+_{G,S}\}$

\item[]
$\check{\Lambda}^+_{M,S}=\{\check{\lambda}\in\check{\Lambda}^+_M\mid
\<w_0^M(\lambda), \check{\lambda}\>\ge 0 \;\, \mbox{for all}\;\,
\lambda\in\Lambda^+_{M,S}\}$ \ \  \QED
\end{itemize}
\end{Lm}

\begin{Def} Let $\Gr_M^+\subset\Gr_M$ be the closed subscheme 
given by the condition: $(\cF_M,\beta)\in\Gr_M^+$ iff
for each $\check{\lambda}\in\check{\Lambda}^+_{M,S}$ the map
$$
\beta^{\check{\lambda}}: U^{\check{\lambda}}_{\cF_M}\to U^{\check{\lambda}}_{\cF^0_M}
$$
is regular on $\cD$. Here $U^{\check{\lambda}}$ stands for the irreducible 
quotient of the Weil module $\cU^{\check{\lambda}}$ for $M$.
\end{Def}
 
 Clearly, $\Gr_M^+$ is a $M(\hat\cO)$-invariant subscheme of $\Gr_M$. For
$\nu\in\Lambda^+_M$ one has the closed subscheme $\Grb^{\nu}_M\subset\Gr_M$ (cf.
\cite{BG}, sect.~3.2). For $\nu\in\Lambda^+_M$ we
have $\Grb_M^{\nu}\subset \Gr^+_M$ iff $\nu\in\Lambda^+_{M,S}$. 

 Recall that for $\mu\in\pi_1(M)$ the connected component
$\Gr_M^{\mu}$ of $\Gr_M$ classifies $(\cF_M,\beta)\in\Gr_M$ such that for
$\check{\lambda}\in \check{\Lambda}_{G,P}$ we have 
$$
\cV^{\check{\lambda}}_{\cF_M^0}(-\<\nu,\check{\lambda}\>)\;\iso\;
\cV^{\check{\lambda}}_{\cF_M}
$$
For $\mu\in\pi_1(M)$ set $\Gr_M^{+,\mu}=\Gr_M^+\cap\Gr_M^{\mu}$. So, $\Gr_M^{+,\mu}$
is nonempty iff $\mu\in\pi_1^+(M)$. It may be shown that $\Gr_M^{+,\mu}$ is connected for each
$\mu\in\pi_1^+(M)$.

 Recall the following definition from (\cite{BG}, sect. 4.3.1). For
$\mu\in\Lambda_{G,P}$ let $S^{\mu}_P\subset\Gr_G$ denote the locally closed subscheme
classifying $(\cF_G,\,\beta:\cF^0_G\mid_{\cD^*}\iso\cF_G\mid_{\cD^*})\in\Gr_G$ such that
the composition
$$
\cL^{\check{\lambda}}_{\cF^0_T}(-\<\mu,\check{\lambda}\>)\to
\cL^{\check{\lambda}}_{\cF^0_T}\to
\cV^{\check{\lambda}}_{\cF^0_G}\toup{\beta}  
\cV^{\check{\lambda}}_{\cF_G}
$$
has niether pole nor zero over $\cD$ for every
$\check{\lambda}\in\check{\Lambda}_{G,P}\cap\check{\Lambda}^+$.

For each $\nu\in\pi_1(G)$ the component
$\Gr_G^{\nu}$ is stratified by
$S^{\mu}_P$ indexed by those $\mu\in\Lambda_{G,P}$ whose image in $\pi_1(G)$ is $\nu$.
Moreover, we have a natural map $t^{\mu}_S: S^{\mu}_P\to
\Gr_M^{\mu}$.

\begin{Lm}
\label{Lm_Gr_M}
For each $\mu\in \Lambda_{G,P}$ the map $t^{\mu}_S: \Gr_G^+\cap S^{\mu}_P\to 
\Gr_M^{\mu}$ factors through $\Gr_M^{+,\mu}\hook{}\Gr_M^{\mu}$, and the induced map
$t^{+,\mu}_S: \Gr_G^+\cap S^{\mu}_P\to \Gr_M^{+,\mu}$ is surjective.
\end{Lm} 
\begin{Prf}
Let $(\cF_G,\,\beta:\cF^0_G\mid_{\cD^*}\to \cF_G\mid_{\cD^*})
\in\Gr_G^+\cap S^{\mu}_P$. So, for any
$\check{\lambda}\in\check{\Lambda}_{G,P}\cap\check{\Lambda}^+$ 
$$
\cL^{\check{\lambda}}_{\cF_T^0}(-\<\mu,\check{\lambda}\>)\hook{} \cV^{\check{\lambda}}_{\cF_G}
$$
is a subbundle. There is unique  $(\cF_P,\beta:\cF_P^0\mid_{\cD^*}\to \cF_P\mid_{\cD^*})\in\Gr_P$ that induces $(\cF_G,\beta)$, and $t^{\mu}_S$ sends $(\cF_G,\beta)$ to $\cF_M=\cF_P\times_P M$. Since for any $\check{\lambda}\in\check{\Lambda}^+_S$ the maps
$$
\beta^{\check{\lambda}}: \cV^{\check{\lambda}}_{\cF_P}\hook{} \cV^{\check{\lambda}}_{\cF^0_P}
$$
are regular, the first assertion is reduced to the next sublemma.

\begin{Slm}
Let $\check{\nu}\in\check{\Lambda}^+_{M,S}$. So, there exists 
$w\in W$ with $w\check{\nu}\in\check{\Lambda}^+_S$. Let $\Res_G^P \,V^{w\check{\nu}}$ denote $V^{w\check{\nu}}$ viewed as a $P$-module. There exists a subquotient $V'$ of 
$\Res_G^P \,V^{\check{\lambda}}$ on which $U(P)$ acts trivially and such that 
$$
\Hom_M(U^{\check{\nu}}, V')\ne 0  \eqno{\square}
$$
\end{Slm}

Let $\nu\in\Lambda^+_{M,S}$ be such that $\Gr_M^{\nu}\subset \Gr_M^{+,\mu}$. 
Recall the notation $\hat\cO=k[[t]]$. Since $t^{+,\mu}_S$ is $M(\hat\cO)$-invariant, it
suffices to show that $\nu(t)\in\Gr_M^{+,\mu}$ lies in the image of $t^{+,\mu}_S$. We
know that there exists $w\in W$ with $w\nu\in\Lambda^+_{G,S}$. Therefore,
$\nu(t)G(\hat\cO)$ defines a point of
$\Gr_G^+\cap S^{\mu}_P$ which is sent by $t^{+,\mu}_S$ to $\nu(t)\in\Gr_M^{+,\mu}$.\\
\end{Prf}(Lemma~\ref{Lm_Gr_M})

\bigskip

 Note as a consequence that for each $\nu\in\pi_1^+(G)$ the scheme $\Gr_G^{+,\nu}$ 
is stratified by locally closed subschemes $\Gr_G^{+,\nu}\cap S^{\mu}_P$ indexed by
those $\mu\in\pi_1^+(M)$ whose image in $\pi_1(G)$ is $\nu$.

 For $d\ge 0$ write $\Lambda^{+,d\theta}_{M,S}$ for the preimage of $d\theta$ under
$\Lambda^+_{M,S}\to\pi_1^+(G)$.

\begin{Lm}
For any $\lambda\in \Lambda^{+,d\theta}_{M,S}$ 
there exist $\lambda_1,\ldots,\lambda_d\in \Lambda^{+,\theta}_{M,S}$ such that 
$\lambda\,\leM\,
\lambda_1+\ldots+\lambda_d$.
\end{Lm}   
\begin{Prf}
Pick any $k$-point $(\cF_M,\beta)$ of $\Gr_M^{\lambda}$. Let $\mu$ be the image of 
$\lambda$ in $\pi_1^+(M)$. Pick any $k$-point $(\cF_G, \beta)$ of $\Gr_G^{+,d}\cap
S^{\mu}_P$ whose $t^{+,\mu}_S$-image is $(\cF_M,\beta)$. Let 
$\Gr^{\gamma}_G\ttimes\ldots\ttimes\Gr_G^{\gamma}$ be the scheme classifying
collections
\begin{equation}
\label{point_1}
(\cF^1_G,\ldots,\cF^{d+1}_G=\cF^0_G,\,\beta_i),
\end{equation}
where $\cF^i_G$ is a $G$-bundle on $\cD$, and $\beta^i:
\cF^i_G\mid_{\cD^*}\,\iso\,\cF^{i+1}_G\mid_{\cD^*}$ is an isomorphism 
such that $\cF^i_G$ is in the position $\gamma$ with respect to $\cF^{i+1}_G$ for
$i=1,\ldots,d$.

 Pick a $k$-point (\ref{point_1}) 
whose image under the convolution map 
$$
\Gr^{\gamma}_G\ttimes\ldots\ttimes\Gr_G^{\gamma}\to\Grb^{d\gamma}_G,
$$ 
is $(\cF_G, \beta)$. There exist a unique collection
$$
(\cF_P^1,\ldots,\cF_P^d, \cF_P^{d+1}=\cF^0_P, \beta_i),
$$
where $\cF^i_P$ is a $P$-torsor on $\cD$ and 
$\beta_i: \cF^i_P\mid_{X-x_i}\iso\,\cF^{i+1}_P\mid_{X-x_i}$, that induces
(\ref{point_1}) by extension of scalars from $P$ to $G$. 
Extending the scalars from $P$ to $M$, one gets a collection
$$
(\cF_M^1,\ldots,\cF_M^d, \cF_M^{d+1}=\cF^0_M, \beta_i)
$$
where $\cF_M^1=\cF_M$. For $i=1,\ldots,d$ let $\lambda_i\in\Lambda^+_M$ be such that
$\cF_M^i$ is in the position $\lambda_i$ with respect to $\cF^{i+1}_M$.
By Lemma~\ref{Lm_Gr_M}, $\lambda_i\in\Lambda^{+,\theta}_{M,S}$ for all $i$.
Let 
$$
\Gr_M^{\lambda_1}\ttimes \ldots\ttimes\Gr_M^{\lambda_d}
$$ 
be the scheme classifying
$(\cF^1_M,\ldots,\cF^{d+1}_M=\cF^0_M, \beta^i)$, where $\cF^i_M$ is a $M$-torsor on
$\cD$ and $\beta^i:\cF^i_M\mid_{\cD^*}\,\iso\,\cF^{i+1}_M\mid_{\cD^*}$ is an isomorphism
such that $\cF^i_M$ is in the position $\lambda_i$ with respect to $\cF^{i+1}_M$ for
$i=1,\ldots,d$.
We learn that $(\cF_M,\beta)$ lies in the image of the convolution map
$$
\Gr_M^{\lambda_1}\ttimes \ldots\ttimes\Gr_M^{\lambda_d}\to \Gr_M^{+,\mu}
$$
But the image of the latter map is contained in $\Grb_M^{\lambda_1+\ldots+\lambda_d}$,
so $\Gr_M^{\lambda}\subset \Grb_M^{\lambda_1+\ldots+\lambda_d}$.
\end{Prf}

\bigskip

Denote by $\pi_1^{\theta}(M)$ the image of $\Lambda^{+,\theta}_{M,S}\to\pi_1^+(M)$.
By the above lemma, $\pi_1^{\theta}(M)$ generates 
$\pi_1^+(M)$ as a semigroup. Since $V^{\gamma}$ is faithful, $\pi_1^{\theta}(M)$
generates $\pi_1(M)$ as a group. 

\begin{Lm} 1) Each $\lambda\in\Lambda^{+,\theta}_{M,S}$ is a minuscule dominant coweight
for $M$. \\
2) The natural map $\Lambda^{+,\theta}_{M,S}\to\pi_1^{\theta}(M)$ is
bijective.
\end{Lm}
\begin{Prf} 1) For $\lambda\in\Lambda^{+,\theta}_{M,S}$ we have
$\Hom_M(U^{\lambda}, V^{\gamma})\ne 0$. Let
$\lambda'\in\Lambda^+_M$ and
$\lambda'\leM\lambda$. Then $\lambda'$ is a weight of $V^{\gamma}$, so that
$\lambda'=w\gamma$ for some $w\in W$, and $\lambda=w\gamma+\alpha$, where $\alpha$ is a sum
of positive coroots for $M$. However, if $\alpha\ne 0$ then the length of $\lambda$ is
strictly bigger than the length of $\gamma$, which implies that $\lambda$ is not a weight
of $V^{\gamma}$. This contradiction shows that $\lambda'=\lambda$.\\
2) follows from 1).
\end{Prf}

\bigskip\noindent
4.2 For $d\ge 0$ consider the stack
\begin{equation}
\label{stack_17}
\cH^{+,d}_G\times_{\Bun_G}\Bun_P,
\end{equation}
where we used the projection $\gq:\cH^{+,d}_G\to\Bun_G$
in the fibred product. 

  For $\mu\in\Lambda_{G,P}$ let $\cH^{+,\mu}_P$ be the locally closed
substack of (\ref{stack_17}) classifying 
$$
(D\in X^{(d)}, \cF_G,\cF'_P,\;\beta:\cF_G\mid_{X-D}\iso\,\cF'_G\mid_{X-D},\;
\cF'_P\times_P G\,\iso\,\cF'_G)
$$ 
for which there exists a $\Lambda_{G,P}$-valued divisor $D^{\mu}$ on $X$ of degree $\mu$
with the property: for all $\check{\lambda}\in\check{\Lambda}_{G,P}\cap\check{\Lambda}^+$ the
meromorphic maps
$$
\cL^{\check{\lambda}}_{\cF'_{M/[M,M]}}\to \cV^{\check{\lambda}}_{\cF'_G}\to
\cV^{\check{\lambda}}_{\cF_G}
$$
realise $\cL^{\check{\lambda}}_{\cF'_{M/[M,M]}}$ as a subundle of 
$\cV^{\check{\lambda}}_{\cF_G}(\<D^{\mu}, \check{\lambda}\>)$. Here we have denoted
$\cF'_M=\cF'_P\times_P M$.

 Clearly, we have $\<D^{\mu}, \check{\omega}_0\>=D$ and
$\<D^{\mu},\check{\omega}_i\>\ge 0$ for $i\in{\cal I}-{\cal I}_M$. Let
$D_i=\<D^{\mu},\check{\omega}_i\>$ for $i\in{\cal I}-{\cal I}_M$ then we have an
equality of $\pi_1(M)$-valued divisors on $X$
$$
D^{\mu}=w_0(\gamma)D+
\mathop{\sum}_{i\in{\cal I}-{\cal I}_M} D_i\alpha_i
$$ 
 
 By Lemma~\ref{Lm_Gr_M}, $\cH^{+,\mu}_P$ is non empty iff $\mu\in\pi_1^+(M)$ and
actually $D^{\mu}$ is a $\pi_1^+(M)$-valued divisor on $X$. So, for each $d\ge 0$ the
stack (\ref{stack_17}) is stratified by locally closed substacks $\cH^{+,\mu}_P$
indexed by those $\mu\in\pi_1^+(M)$ whose image in $\pi_1(G)$ is $d\theta$.

  For $\mu\in\pi_1^+(M)$ let $d=\<\mu,\check{\omega}_0\>$ and
$d_i=\<\mu,\check{\omega}_i\>$ for
$i\in{\cal I}-{\cal I}_M$ and let $X^{\mu}_M$ denote the scheme image of the projection
$$
\cH^{+,\mu}_P\to X^{(d)}\times\prod_{i\in{\cal I}-{\cal
I}_M} X^{(d_i)}
$$
We will think of $X^{\mu}_M$ as the moduli scheme of $\pi_1^+(M)$-valued 
divisors on $X$ of degree $\mu$.  
As we will see, $X^{\mu}_M$ need not be irreducible. 
For $\mu\in\pi_1^+(M)$ we have a commutative diagram 
$$
\begin{array}{ccc}
\cH^{+,\mu}_P & \;\;\hook{}\; & \cH^{+,d}_G\times_{\Bun_G}\Bun_P\\
\downarrow\lefteqn{\scriptstyle\supp_P} && \downarrow\\
X^{\mu}_M & \;\;\toup{s_M}\; & X^{(d)},
\end{array}
$$
where we have denoted by $\supp_P$ and $s_M$ the natural projections.

\smallskip

\begin{Def} For $\mu\in\pi_1^+(M)$ let $\cH^{+,\mu}_M$ be the stack of collections
$(\cF_M,\cF'_M, D^{\mu}\in X^{\mu}_M, \beta)$, where $\cF_M, \cF'_M\in\Bun_M$ and
$$
\beta: \cF_M\mid_{X-D}\iso \cF'_M\mid_{X-D}
$$
is an isomorphism of $M$-torsors with
$D=s_M(D^{\mu})$ such that for each $\check{\lambda}\in\check{\Lambda}^+_{M,S}$ the map
$$
\beta^{\check{\lambda}}: U^{\check{\lambda}}_{\cF_M}\hook{} U^{\check{\lambda}}_{\cF'_M}
$$
extends to an inclusion of coherent sheaves on $X$, and 
$$
\cL^{\check{\lambda}}_{\cF_{M/[M,M]}}(\<D^{\mu},\check{\lambda}\>)\,\iso\, 
\cL^{\check{\lambda}}_{\cF'_{M/[M,M]}}
$$
for each $\check{\lambda}\in\check{\Lambda}_{G,P}$.
\end{Def}

We have a diagram 
$$
\begin{array}{ccc}
\Bun_M\,\getsup{\gp_M} & \cH^{+,\mu}_M & \toup{\gq_M}\,\Bun_M\\
& \downarrow\lefteqn{\scriptstyle \,\supp_M}\\
& X^{\mu}_M
\end{array}
$$
where $\gp_M$ (resp., $\gq_M$) sends a point of $\cH^{+,\mu}_M$ to $\cF_M$ 
(resp., to $\cF'_M$), and $\supp_M$ stands for the projection.
If $(\cF'_M, D^{\mu}=\sum_k \mu_k x_k)$ is a $k$-point of
$\Bun_M\times X^{\mu}_M$ then the fibre of 
$$
\gq_M\times\supp_M:\cH^{+,\mu}_M\to \Bun_M\times X^{\mu}_M
$$ 
over it identifies with
$\prod_k \Gr_M^{+,\mu_k}$.

\bigskip\noindent
4.3 Given $\mu\in\pi_1^+(M)$, we will denote by $\gA(\mu)$ the elements of the set 
of decompositions of $\mu$ as a sum of non-zero elements of $\pi_1^+(M)$. More
precisely, $\gA(\mu)$ is a way to write $\mu=\sum_k n_k\mu_k$, where all $n_k>0$ and 
$\mu_k\in\pi_1^+(M)-\{0\}$ are pairwise distinct.

 For $\gA(\mu)$ we set $X^{\gA(\mu)}=\prod_k X^{(n_k)}$. We have a natural map 
$s^{\gA(\mu)}: X^{\gA(\mu)}\to X^{\mu}_M$ sending $\{D_k\}$ to $\sum_k \mu_kD_k$.  Let
$\Xo^{\gA(\mu)}\subset X^{\gA(\mu)}$ be the complement to all diagonals.  
The composition 
$$
\Xo^{\gA(\mu)}\hook{} X^{\gA(\mu)}\to X^{\mu}_M
$$
is a locally closed embedding, and in this way $X^{\mu}_M$ is stratified by subschemes
$\Xo^{\gA(\mu)}$.

 We say that $\gA(\mu)$ \select{is in general position} if  $\sum_k n_k=d$.
Write $^{rss}X^{\mu}_M$ for the preimage of $^{rss}X^{(d)}$ under 
$s_M: X^{\mu}_M\to X^{(d)}$. The connected components of $^{rss}X^{\mu}_M$ are exactly 
$\Xo^{\gA(\mu)}$ indexed by $\gA(\mu)$ in general position.

\begin{Def} Given a local system $W$ on $X$, for each $\mu\in\pi_1^+(M)$ define 
Laumon's sheaf
$\cL^{\mu}_W$ on $\cH^{+,\mu}_M$ as follows. Recall the diagram
$$
\cH^{+,\mu}_M\toup{\supp_M} X^{\mu}_M\toup{s_M} X^{(d)}
$$
Let $^{rss}\cH^{+,\mu}_M$ be the preimage of $^{rss}X^{(d)}$ under $s_M\comp\supp_M$.
The stack $^{rss}\cH^{+,\mu}_M$ is smooth and over it
we let
$$
\cL^{\mu}_W=\supp_M^*s_M^*W^{(d)}[a](\frac{a}{2}),
$$
where
$a$ denotes the dimension of the corresponding connected component of 
$^{rss}\cH^{+,\mu}_M$. 
Then we extend this perverse sheaf by Goresky-MacPherson to $\cH^{+,\mu}_M$. 
\end{Def}

By definition,
$\DD(\cL^{\mu}_W)\iso\cL^{\mu}_{W^*}$, and
$\cL^{\mu}_W\;\iso\;\oplus \cL^{\gA(\mu)}_W$ is a direct sum of perverse
sheaves indexed by $\gA(\mu)$ in general position.
Set $d_M=\dim\Bun_M$. 
 
 For  $\mu\in\pi_1^+(M)$ whose image in $\pi_1(G)$ is $d\theta$ 
consider the diagram
$$
\cH^{+,\mu}_M\times_{\Bun_M} \Bun_P \;\getsup{q_M}  \;\cH^{+,\mu}_P  \; \toup{f_M}  
\; \cH^{+,d}_G
$$
where we used $\gq_M:\cH^{+,\mu}_M\to\Bun_M$ in the fibred product, $f_M$ 
is the natural map, and $q_M$ sends $(\cF_P,\cF'_P, D^{\mu},\beta)$ to $(\cF_M, \cF'_P,
D^{\mu}, \; \cF'_P\times_P M\,\iso\,\cF'_M, \; \beta)$. We also have
$$
\cH^{+,\mu}_M\times_{\Bun_M} \Bun_P \;\getsup{p_M}  \;\cH^{+,\mu}_P  \; \toup{f_M}  
\; \cH^{+,d}_G,
$$
where now we used $\gp_M: \cH^{+,\mu}_M\to\Bun_M$ in the fibred product, and
$p_M$ sends $(\cF_P,\cF'_P, D^{\mu},\beta)$ to $(\cF'_M, \cF_P, D^{\mu}, \; 
\cF_P\times_P M\,\iso\,\cF_M, \; \beta)$.

Here is a generalization of Laumon's theorem (\cite{La2}, Theorem~4.1).

\begin{Pp} 
\label{Pp_Laumon_generalized}
Let $W$ be a local system on $X$. 
Let $\mu\in\pi_1^+(M)$ with image $d\theta$ in
$\pi_1^+(G)$.  
The complex $q_{M!}f_M^*\cL^d_W$ is canonically isomorphic to
the inverse image of 
$$
\cL^{\mu}_W\otimes\Qlb[1](\frac{1}{2})^{\otimes d_G-d_M+\<\mu, \,2\check{\rho}_M-2\check{\rho}\>}
$$ 
under the projection
$\cH^{+,\mu}_M\times_{\Bun_M}\Bun_P\to\cH^{+,\mu}_M$. (We have used the fact that
$\check{\rho}_M-\check{\rho}\in\check{\Lambda}_{G,P}$). The complex 
$p_{M!}f_M^*\cL^d_W$ is canonically isomorphic to
the inverse image of 
$$
\cL^{\mu}_W\otimes\Qlb[1](\frac{1}{2})^{\otimes d_G-d_M-\<\mu,
\,2\check{\rho}_M-2\check{\rho}\>}
$$ 
under the projection $\cH^{+,\mu}_M\times_{\Bun_M}\Bun_P\to\cH^{+,\mu}_M$.
\end{Pp}

 The proof is given in Sections~4.4-4.5.

\bigskip\noindent
4.4 Let $J=\{i\in{\cal I}\mid \<\gamma, \check{\alpha}_i\> =0\}$. Let $W_J\subset W$ be
the subgroup generated by the reflection corresponding to $i\in J$.
Using Bruhat decomposition, one checks that the map $W/W_J\to W\gamma$ sending $w$ to
$w\gamma$ is a bijection. 

 Fix a section $T\to B$. Let $P_{\gamma}$ denote the parabolic
of $G$  generated by $T$ and
$U_{\check{\alpha}}$ for all roots $\check{\alpha}$ such that
$\<\gamma,\check{\alpha}\>\le 0$.  So, $P_{\gamma}$ contains the opposite Borel.
We have a bijection $\Lambda^{+,\theta}_{M,S}\,\iso\,
W_M\backslash W/W_J$ sending 
$w\gamma\in \Lambda^{+,\theta}_{M,S}$ to the coset $W_MwW_J$.

  The map $G/P_{\gamma}\to \Gr_G^{\gamma}$ sending $g\in G(k)\subset G(\hat\cO)$ to
$g\gamma(t)G(\hat\cO)$ is an isomorphism. 
The scheme $\Gr_G^{\gamma}$ is stratified by $\Gr_G^{\gamma}\cap S^{\lambda}_P$  
indexed by $\lambda\in\Lambda^{+,\theta}_{M,S}$. The above isomorphism transforms this
stratification into the stratification of $G/P_{\gamma}$ by $P$-orbits. We have a
disjoint decomposition
$$
G=\mathop{\sqcup}\limits_{w\in W_M\backslash W/W_J} P wP_{\gamma}
$$
So,  we have
$$
\Gr_G^{\gamma}\cap S^{\lambda}_P\;\iso\; PwP_{\gamma}/P_{\gamma}\;\iso\;
P/P\cap wP_{\gamma}w^{-1}
$$
 
 Similarly, for $\lambda\in \Lambda^{+,\theta}_{M,S}$ let $P_{\lambda}(M)$ be the 
parabolic of $M$ generated by $T$ and $U_{\check{\alpha}}$, where 
$\check{\alpha}$ runs through those roots of $M$ for which $\<\lambda, 
\check{\alpha}\>\le 0$. Then the map $M/P_{\lambda}(M)\to \Gr_M^{\lambda}$ sending
$m\in M(k)$ to $m\lambda(t)M(\hat\cO)$ is an isomorphism. So, the map
\begin{equation}
\label{map_111}
\Gr_G^{\gamma}\cap S^{w\gamma}_P\to \Gr_M^{w\gamma}
\end{equation}
is nothing else but the map $P/P\cap wP_{\gamma}w^{-1}\to M/P_{w\gamma}(M)$
sending $p$ to $p\!\mod\! M$. The correctness is due to

\begin{Lm} We have $P_{w\gamma}(M)=M\cap wP_{\gamma}w^{-1}$.
\end{Lm}
\begin{Prf} The inclusion $P_{w\gamma}(M)\subset M\cap wP_{\gamma}w^{-1}$
follows from definitions. Further, $M\cap wP_{\gamma}w^{-1}$ contains the opposite 
Borel of $M$, hence is a parabolic subgroup of $M$ (in particular, connected). The
assertion follows now from: for a root $\check{\alpha}$ of $M$ we have
$U_{\check{\alpha}}\subset P_{w\gamma}(M)$ if and only if $U_{\check{\alpha}}\subset
M\cap wP_{\gamma}w^{-1}$. 
\end{Prf}

\medskip\smallskip

 We see that $U(P)$ acts transitively on the fibres of (\ref{map_111}). So, 
(\ref{map_111}) is a fibration with fibre isomorphic 
to an affine space of dimension $\<\gamma+w\gamma,\check{\rho}\>-\<w\gamma, 
2\check{\rho}_M\>$.

\begin{Rem} For $w\gamma\in\Lambda^+_M$
one can calculate the dimension of $\Gr_G^{\gamma}\cap S^{w\gamma}_P$ as follows. 
Stratify it by $B$-orbits, that is, by the schemes
$\Gr_G^{\gamma}\cap S^{w_1w\gamma}$ with $w_1\in W_M$. Then $\Gr_G^{\gamma}\cap 
S^{w_1w\gamma}$ is an affine space of dimension $\<\gamma+w_1w\gamma,\check{\rho}\>$.
The maximum of these numbers, as $w_1$ ranges through $W_M$ is
$\<\gamma+w\gamma,\check{\rho}\>$. 
\end{Rem}

\bigskip\noindent
4.5 Consider a collection $\tilde\mu=(\mu_1,\ldots,\mu_d)$ with
$\mu=\mu_1+\ldots+\mu_d$ and $\mu_i\in\pi_1^{\theta}(M)$. Let
$\cH^{+,\tilde\mu}_M$ be the stack of collections 
\begin{equation}
\label{collection_18}
(\cF^1_M,\ldots,\cF^{d+1}_M,
x_1,\ldots,x_d\in X, \beta^i),
\end{equation}
where $\beta^i: \cF^i_M\mid_{X-x_i}\iso\,
\cF^{i+1}_M\mid_{X-x_i}$ is an isomorphism such that $(\cF^i_M,\cF^{i+1}_M,\beta^i,
x_i)\in\cH^{+,\mu_i}_M$ for
$i=1,\ldots,d$. 

 If we denote by $\lambda_i$ the element of $\Lambda^{+,\theta}_{M,S}$ that maps
to $\mu_i$ then $\cF^i_M$ is in the position $\lambda_i$ with respect to $\cF^{i+1}_M$
at $x_i$.  We have a convolution map 
$$
\conv^{\tilde\mu}: \cH^{+,\tilde\mu}_M\to
\cH^{+,\mu}_M
$$ 
sending (\ref{collection_18}) to $(\cF^1_M,\cF^{d+1}_M, \beta,
D^{\mu})$, where $D^{\mu}=\sum_i \mu_i x_i$ and
$\beta:\cF^1_M\mid_{X-D}\,\iso\,\cF^{d+1}_M\mid_{X-D}$ with $D=s_M(D^{\mu})$. 

 Denote by $s^{\tilde\mu}:
\cH^{+,\tilde\mu}_M\to X^d$ the map sending (\ref{collection_18}) to
$(x_1,\ldots,x_d)$. From
(Lemma~9.3, \cite{NP}) one derives

\begin{Lm} 
\label{Lm_small_map_again}
i) The map $\conv^{\tilde\mu}$ is representable, proper and small over its image. 
Besides, the perverse sheaf 
\begin{equation}
\label{sheaf_1}
\conv^{\tilde\mu}_!(s^{\tilde\mu})^*W^{\boxtimes d}[a](\frac{a}{2})
\end{equation}
is the Goresky-MacPherson extension from $^{rss}\cH^{+,\mu}_M$. Here $a=\dim
\cH^{+,\tilde\mu}_M$. \\
ii) The $d$-tuple $\tilde\mu$ gives rise to $\gA(\mu)$ in general position, say $\mu=\sum_k n_k\nu_k$. The group $\prod_k S_{n_k}$ acts naturally on (\ref{sheaf_1}), and 
the sheaf of  $\prod_k S_{n_k}$-invariants is canonically isomorphic to the direct summand of $\cL^{\mu}_W$ corresponding to $\gA(\mu)$.
\QED
\end{Lm}

\smallskip

\begin{Rem} 
\label{remark_ULA}
From Lemma~\ref{Lm_small_map_again} it follows that for any
$\mu\in\pi_1^+(M)$ the complex $\cL^{\mu}_W$ is ULA with respect to both projections
$\gp_M, \gq_M: \cH^{+,\mu}_M\to\Bun_M$.
\end{Rem}

\medskip
\begin{Prf}\select{of Proposition~\ref{Pp_Laumon_generalized}}\\
1) Consider the diagram 
$$
\begin{array}{ccc}
\cH^{+,\mu}_P & \toup{q_M} & \;\cH^{+,\mu}_M\times_{\Bun_M}\Bun_P\\
& \searrow\lefteqn{\scriptstyle \gq_P\times \supp_P} & \;\downarrow\\
&& \;\Bun_P\times X^{\mu}_M
\end{array}
$$
For any $k$-point of $\Bun_P\times X^{\mu}_M$ given by $(\cF'_P,\; D^{\mu}=\sum_k \mu_k
x_k)$, where $x_k$ are pairwise distinct,
over $(\cF'_P, D^{\mu}=\sum_k \mu_k x_k)$ the map $q_M$ becomes the product
$$
\prod_k \Gr_G^+\cap S^{\mu_k}_P\to \prod_k \Gr_M^{+,\mu_k}
$$
of the maps $t^{+,\mu_k}_S$. In particular, $q_M$ is surjective.

 Let $\gA(\mu)$ in general position be given by $\mu=\sum_k n_k\mu_k$. Let 
$U\subset {^{rss}\cH^{+,\mu}_M}$ be
the preimage of the corresponding connected component of $^{rss}X^{\mu}_M$.
Over the open substack $U\times_{\Bun_M}\Bun_P$, the map $q_M$ is a
fibration with fibre isomorphic to an affine space of dimension
$$
a(\gA(\mu))\df\<d\gamma,\check{\rho}\>+
\sum_k \< n_k\lambda_k, \,\check{\rho}-2\check{\rho}_M\>,
$$
where $\lambda_k\in\Lambda^{+,\theta}_{M,S}$ maps to $\mu_k\in\pi_1^{\theta}(M)$.

 The restriction of $f_M^*\cL^d_W$ to $q_M^*(U\times_{\Bun_M}\Bun_P)$ comes from 
$^{rss}X^{(d)}$. So, over $U\times_{\Bun_M}\Bun_P$, we get the desired isomorphism.  
Now it suffices to show that, up to a shift, $q_{M!}f_M^*\Spr^d_W$ is a perverse sheaf, 
the Goresky-MacPherson extension from
$^{rss}\cH^{+,\mu}_M\times_{\Bun_M}\Bun_P$. 

 For a $d$-tuple $\tilde\mu=(\mu_1,\ldots,\mu_d)$ with  $\mu=\mu_1+\ldots+\mu_d$  and 
$\mu_i\in\pi_1^{\theta}(M)$, let $\cH^{+,\tilde\mu}_P$ be the stack of collections
\begin{equation}
\label{tuple_1}
(\cF^1_P,\ldots,\cF^{d+1}_P, x_1,\ldots,x_d,\beta^i),
\end{equation}
where $x_i\in X$ and
$(\cF^i_P,\cF^{i+1}_P, x_i,\beta^i)\in\cH^{+,\mu_i}_P$ for $i=1,\ldots,d$.
The stack 
$$
\cH^{+,\mu}_P\times_{(\cH^{+,d}_G\times_{\Bun_G}\Bun_P)} 
(\tilde\cH^{+,d}_G\times_{\Bun_G}\Bun_P)
$$ 
is stratified by locally closed substacks
$\cH^{+,\tilde\mu}_P$ indexed by such tuples $\tilde\mu$.

We have a diagram 
$$
X^{(d)}\,\getsup{s^{\tilde\mu}}\,\cH^{+,\tilde\mu}_P\toup{\alpha}
\cH^{+,\tilde\mu}_M\times_{\Bun_M}\Bun_P\to
\cH^{+,\mu}_M\times_{\Bun_M}\Bun_P,
$$
where $\alpha$ sends (\ref{tuple_1}) to $(\cF^1_M,\ldots, \cF^d_M, \cF^{d+1}_P,
x_1,\ldots,x_d, \beta^i)$ with $\cF^i_M=\cF^i_P\times_P M$, and the last map is
$\conv^{\tilde\mu}\times\id$. 
It suffices to show that for each $\tilde\mu$ as above,
$$
(\conv^{\tilde\mu}\times\id)_!\alpha_!(s^{\tilde\mu})^*W^{(d)}
$$
is a perverse sheaf (up to a shift), the Goresky-MacPherson
extension from $^{rss}\cH^{+,\mu}_M\times_{\Bun_M}\Bun_P$.
Since $\alpha$ is a composition of affine fibrations, our statement about
$q_{M!}f_M^*\cL^d_W$ follows from
Lemma~\ref{Lm_small_map_again}, i).

\smallskip\noindent
2) Applying 1) to the 1-admissible data $\{-w_0(\gamma)\}$, one gets the formula for
$p_{M!}f_M^*\cL^d_W$.
\end{Prf}

\bigskip\noindent
4.6.1 \select{Averaging functors for Levi subgroups} \\
Recall for any $\mu\in\pi_1^+(M)$ the
diagram 
$\Bun_M\,\getsup{\gp_M} \cH^{+,\mu}_M \toup{\gq_M}\,\Bun_M$.
For a local system $W$ on $X$ denote by $\Av^{\mu}_W: \D(\Bun_M)\to\D(\Bun_M)$ the
functor
$$
\Av^{\mu}_W(K)=(\gq_M)_!(\gp_M^*K\otimes\cL^{\mu}_W)[-d_M](\frac{-d_M}{2})
$$
Let also $\Av^{-\mu}_W: \D(\Bun_M)\to\D(\Bun_M)$ be given by 
$$
\Av^{-\mu}_W(K)=(\gp_M)_!(\gq_M^*K\otimes\cL^{\mu}_W)[-d_M](\frac{-d_M}{2})
$$
By Remark~\ref{remark_ULA}, we have $\DD\comp \Av^{\mu}_W\iso \Av^{\mu}_{W^*}\comp\,
\DD$ and
$\DD\comp \Av^{-\mu}_W\iso \Av^{-\mu}_{W^*}\comp\, \DD$ naturally. 

 The proof of the following result is completely analogous to that of
Proposition~\ref{Pp2}.

\begin{Pp} 
\label{Pp_averaging}
Let $\mu\in \pi_1^+(M)$. Let $\gA(\mu)$ in general position be given by
$\mu=\sum_k n_k\mu_k$. Recall the map $s^{\gA(\mu)}: X^{\gA(\mu)}\to X^{\mu}_M$ (cf.
sect.~4.3).  Let $\lambda_k\in \Lambda^{+,\theta}_{M,S}$ be the
element that maps to $\mu_k\in\pi_1^{\theta}(M)$.  
Let $W$ be any local system on $X$. Let $E$ be a $\check{M}$-local system on $X$, 
$K$ be $E$-Hecke eigensheaf on $\Bun_M$. Then for the diagram
$$
\Bun_M\getsup{\gp_M}\cH^{+,\mu}_M\toup{\gq_M\times\supp_M} \Bun_M\times X^{\mu}_M
$$
there is a canonical isomorphism
$$
(\gq_M\times\supp_M)_!(\gp_M^*K\otimes \cL^{\gA(\mu)}_W)[-d_M](\frac{-d_M}{2})\;\iso\;
K\boxtimes s^{\gA(\mu)}_!(\mathop{\boxtimes}\limits_k (W\otimes U^{\lambda_k}_E)^{(n_k)})
[d](\frac{d}{2})
\eqno{\square}
$$
\end{Pp}

\begin{Cor} For each standard proper parabolic $P$ of $G$ there exists a constant $c(P)$
with the following property. Let $E$ be any
$\check{M}$-local system $E$ on $X$, $K$ be a $E$-Hecke eigensheaf on $\Bun_M$. 
Let $W$ be an irreducible local system on $X$ of rank $r=\dim V^{\gamma}$.
For any $\mu\in\pi_1^+(M)$ whose image in $\pi_1(G)$ is $d\theta$ with $d>c(P)$ 
we have $\Av^{\mu}_W(K)=0$.
\end{Cor}
\begin{Prf}
The functor $\Av^{\mu}_W=\oplus \Av^{\gA(\mu)}_W$ is a direct sum of functors indexed by
$\gA(\mu)$ in general position. In the notation of Proposition~\ref{Pp_averaging},
we have
\begin{equation}
\label{iso_avegaring_acts}
\Av^{\gA(\mu)}_W(K)\,\iso\, K\otimes (\mathop{\otimes}\limits_k \Sym^{n_k}
\RG(X, W\otimes U^{\lambda_k}_E))[d](\frac{d}{2})
\end{equation}
Here $d=\sum n_k$, and $k$ runs through the finite set $\pi_1^{\theta}(M)$. For $d$
large enough at least one of $n_k$ will satisfy $n_k>r(2g-2)\dim U^{\lambda_k}$, and the
RHS of (\ref{iso_avegaring_acts}) will vanish.
\end{Prf}

\bigskip

 Generalizing the Vanishing Conjecture of Frenkel, Gaitsgory and Vilonen (\cite{FGV2}),
we suggest
\begin{Con}
\label{Con_3}
Let $W$ be an irreducible local system on $X$ of rank $r=\dim V^{\gamma}$. 
Assume that $P$ is a standard proper parabolic of $G$. 
Then for all $\mu\in\pi_1^+(M)$ whose image in $\pi_1(G)$
equals $d\theta$ with $d>c(P)$, the functor $\Av^{\mu}_W$ vanishes identically. 
\end{Con}

\noindent
4.6.2  Consider the diagram $\Bun_G \getsup{\alpha_P}\Bun_P\toup{\beta_P}\Bun_M$, where
$\alpha_P$ and $\beta_P$ are natural maps. The constant term functor
$\CT_P:\D(\Bun_G)\to \D(\Bun_M)$ is defined by $\CT_P(K)=\beta_{P!}\alpha^*_P(K)$.

 The following is a generalization of Lemma~9.8,\cite{FGV2}.

\begin{Lm} 
\label{Lm_CT_filtration}
Let $W$ be any local system on $X$. For any $K\in\D(\Bun_G)$ and $d\ge 0$ the
complex
$\CT_P\comp\Av^d_W(K)\in\D(\Bun_M)$ has a canonical filtration by complexes
\begin{equation}
\label{contribution_1}
\Av^{\mu}_W\comp\CT_P(K)\otimes\Qlb[1](\frac{1}{2})^{\otimes \<\mu,
\,2\check{\rho}-2\check{\rho}_M\>}
\end{equation}
indexed by those $\mu\in\pi_1^+(M)$ whose image in $\pi_1(G)$ is $d\theta$.
\end{Lm}
\begin{Prf} 
Consider the stack $\cH^{+,d}_G\times_{\Bun_G}\Bun_P$, where we used $\gq:\cH^{+,d}_G\to
\Bun_G$ in the fibred product. The complex $\CT_P\comp\Av^d_W(K)$ is the direct image
with respect to the natural map
\begin{equation}
\label{map_CT}
\cH^{+,d}_G\times_{\Bun_G}\Bun_P\to\Bun_M
\end{equation}
Recall that $\cH^{+,d}_G\times_{\Bun_G}\Bun_P$ is stratified by locally closed substacks
$\cH^{+,\mu}_P$ indexed by those $\mu\in\pi_1^+(M)$ whose image in $\pi_1(G)$ is
$d\theta$. This gives a filtration on $\CT_P\comp\Av^d_W(K)$. 

 The restriction of the map (\ref{map_CT}) to the strutum $\cH^{+,\mu}_P$ can be written
as a composition
$$
\cH^{+,\mu}_P \toup{p_M} \cH^{+,\mu}_M\times_{\Bun_M}\Bun_P\to \Bun_M
$$
So, by Proposition~\ref{Pp_Laumon_generalized}, the contribution of the stratum
$\cH^{+,\mu}_P$ to the direct image in question is exactly (\ref{contribution_1}).
\end{Prf}

\begin{Cor} 
\label{Cor4}
Assume that Conjecture~\ref{Con_3} holds. Then \\
1) Let $d$ satisfy $d>c(P)$ for any standard
proper parabolic of $G$. Then for any $K\in\D(\Bun_G)$ and any irreducible
local system $W$ on $X$ of rank $r=\dim V^{\gamma}$ the complex $\Av^d_W(K)$ is
cuspidal. 

\smallskip\noindent
2) Let $E$ be $\check{G}$-local system on $X$ and $K$ be a $E$-Hecke eigensheaf on
$\Bun_G$. If $V^{\gamma}_E$ is irreducible then $K$ is cuspidal.
\end{Cor}
\begin{Prf}
1) is clear.\\
2) The argument given in (\cite{FGV2}, Theorem 9.2) applies in our setting. Namely,
pick 
$d$ such that $d>c(P)$ for any standard
proper parabolic of $G$. Set $W=(V^{\gamma}_E)^*$. 
By Proposition~\ref{Pp2}, 
$$
\CT_P\comp\Av^d_W(K)\,\iso\, \CT_P(K)\otimes\RG(X^{(d)}, (W\otimes V^{\gamma}_E)^{(d)})
[d](\frac{d}{2})
$$
The LHS vanishes by Lemma~\ref{Lm_CT_filtration}. Since $\RG(X^{(d)}, (W\otimes
V^{\gamma}_E)^{(d)})$ is not zero, $\CT_P(K)=0$.
\end{Prf}

\begin{Rem} For $G=\GL_n$ Conjecture~\ref{Con_3} is proved by D. Gaitsgory (\cite{Ga}).
For $G=\GSp_4$ (example 2 in the appendix) Conjecture~\ref{Con_3} also holds, it is
easily reduced to the result of \select{loc.cit.} for $\GL_2$. So, for $G=\GSp_4$ 
Corolary~\ref{Cor4} is unconditional. 
\end{Rem}

\appendix
\bigskip\medskip
\centerline{\scshape Appendix. 1-Admissible groups}      

\begin{Def} Let $H$ be a connected, semi-simple and simply-connected group (over $k$). Assume
that the center $Z(H)$ is cyclic of order $h$ and fix an isomorphism ${\mathfrak z}:\mu_h\iso
Z(H)$. Assume that the characteristic of $k$ does not divide $h$.
Denote by $G$ the quotient of $H\times\Gm$ by the diagonally embedded $\mu_h$.  Call a
reductive group $G$ over $k$ \select{1-admissible}, if it is obtained in this way.
\end{Def}

Let $H$ be a connected, semi-simple and simply-connected group (over $k$).  Let $T_H$ be a
maximal torus of $H$. Write $\check{\Lambda}_H$ (resp.,
${\Lambda}_H$) for the weight (resp., coweight) lattice of $T_H$. Let $\check{Q}_H\subset
\check{\Lambda}_H$ be the root lattice. Set 
$$
Q_H=\{\lambda\in {\Lambda}_H\otimes\QQ\mid \<\lambda,\check{\lambda}\>\in\ZZ\;
\mbox{for any}\; \check{\lambda}\in\check{Q}_H\}
$$ 
We have a natural pairing 
$\check{\Lambda}_H/\check{Q}_H\times Q_H/\Lambda_H\to
(\frac{1}{h}\ZZ)/\ZZ$. Therefore, any isomorphism 
$$
\tau:Q_H/\Lambda_H\,\iso\,
(\frac{1}{h}\ZZ)/\ZZ
$$ 
yields an isomorphism $\check{\tau}:
\check{\Lambda}_H/\check{Q}_H\,\iso\,\ZZ/h\ZZ$. Since
the characteristic of $k$ does not divide $h$, 
$\Hom(Z(H), k^*)\,\iso\, \check{\Lambda}_H/\check{Q}_H$
canonically. So, $\check{\tau}$ yields ${\mathfrak z}:\mu_h\iso Z(H)$ such that for 
$x\in\mu_h, \check{\lambda}\in\check{\Lambda}_H$ we have
$$
\check{\lambda}({\mathfrak z}(x))=x^{\check{\tau}(\check{\lambda})}
$$ 
So, $\tau$ gives rise to a 1-admissible group $G=(H\times\Gm)/\mu_h$. 

 For the maximal torus $T=(T_H\times\Gm)/\mu_h$ in $G$ the weight lattice is
$$
\check{\Lambda}=\{(\check{\lambda},a)\in\check{\Lambda}_H\times\ZZ\mid
\check{\tau}(\check{\lambda})+a=0\!\mod h\}
$$
and the coweight lattice is
$$
\Lambda=\{(\lambda,b)\in Q_H\times (\frac{1}{h}\ZZ)\mid \tau(\lambda)-b\in\ZZ\}
$$
It is understood that the pairing $\Lambda\times\check{\Lambda}\to\ZZ$ sends
$(\lambda,b),(\check{\lambda},a)$ to $\<\lambda,\check{\lambda}\>+ab$. The map
$(\lambda,b)\mapsto b$ yields an isomorphism
$\pi_1(G)\,\iso\, \Lambda/\Lambda_H\,\iso\, \frac{1}{h}\ZZ$. Note also that 
$\pi_1(\check{G})\,\iso\, \check{\Lambda}/\check{Q}_H\,\iso\,\ZZ$.

 The next result follows from definitions.

\begin{Lm} 
\label{Lm_1-adm_examples}
Let $\gamma_H\in Q_H$ be a dominant coweight for $H_{ad}=H/Z(H)$. Assume that 
\begin{itemize}
\item either $\gamma_H$ is minuscule or ($H=1$ and $\gamma_H=0$);
\item  $\gamma_H$ generates $Q_H/\Lambda_H$;
\item the irreducible representation $V^{\gamma_H}$ of
$(H_{ad})^{\check{}}$ is faithful.
\end{itemize}
Let $\tau: Q_H/\Lambda_H\,\iso\,
(\frac{1}{h}\ZZ)/\ZZ$ be the isomorphism sending $\gamma_H$ to $\frac{1}{h}$, and $G$ be the
corresponding 1-admissible group. Set $\gamma=(\gamma_H, \frac{1}{h})\in\Lambda$. Then
$\{\gamma\}$ is a 1-admissible datum for $G$. \QED
\end{Lm}

\bigskip\smallskip\noindent
\centerline{{\scshape Examples of 1-admissible data}}

\bigskip\noindent
The examples below are produced using Lemma~\ref{Lm_1-adm_examples}. 

\medskip
\noindent
1. {\bf The case $G=\GL_n$.} In the standard notation $\Lambda=\ZZ^n$,
$\check{\Lambda}=\ZZ^n$.  For $1\le i<n$ take $\check{\omega}_i=(1,\ldots,1,0,\ldots,0)$ 
where 1 appears $i$ times, and $\check{\omega}_0=(1,\ldots,1)$. All the
conditions are satisfied and $\gamma=(1,0,\ldots,0)$. 

 Let $\gamma_i=(1,\ldots,1,0,\ldots,0)\in\Lambda$, where 1 appears $i$ times. Then 
$\Lambda^+_{G,S}$ is the $\ZZ_+$-span of $\gamma_1,\gamma_2,\ldots,\gamma_n$. So,
$\Lambda^+_{G,S}\,\iso\, (\ZZ_+)^n$. The element $\omega=\gamma_n$ generates the group of
coweights orthogonal to all roots.

\begin{Rem} This particular choice of $\gamma_H$ yields a construction of an automorphic 
sheaf proposed by Laumon in \cite{La}. However, all fundamental coweights for
$H_{ad}=\PSL_n$ are minuscule, and a choice of $\gamma_H$ here is equivalent to a choice
of a generator of the cyclic group $\pi_1(H_{ad})$. If $\gamma_H$ is a fundamental
coweight corresponding to a simple root which is not one of two edges of the Dynkin
diagram $A_{n-1}$ then the corresponding 1-admissible group $G$ is not isomorphic to
$\GL_n$.
\end{Rem}

\medskip\noindent
2. {\bf The case $G=G\Sp_{2n}$, $n\ge 1$.} The group $G$ is a quotient of
$\Gm\times\Sp_{2n}$ by the diagonally embedded $\{\pm 1\}$. Realise $G$ as the subgroup of
$\GL(k^{2n})$ preserving up to a scalar the bilinear form given by the matrix
$$
\left(
\begin{array}{cc}
0 & E_n\\
-E_n & 0
\end{array}
\right),
$$
where $E_n$ is the unit matrix of $\GL_n$.

 The maximal torus $T$ of $G$ is $\{(y_1,\ldots,y_{2n})\mid y_iy_{n+i} 
\mbox{ does not depend on  } i\}$. Let $\check{\epsilon}_i\in\check{\Lambda}$ 
be the caracter that sends a point of $T$ to $y_i$. 
The roots are
$$
\check{R}=\{\pm\check{\alpha}_{ij} \; (i< j \in 1,\ldots,n),\;
\pm\check{\beta}_{ij}\; (i\le j \in 1,\ldots,n)\},
$$
where
$\check{\alpha}_{ij}=\check{\epsilon}_i-\check{\epsilon}_j$ and 
$\check{\beta}_{ij}=\check{\epsilon}_i-\check{\epsilon}_{n+j}$. 

 We have $\Lambda=\{(a_1,\ldots,a_{2n})\mid a_i+a_{n+i} \mbox{ does not depend on }
i\}$. The weight latice is 
$$
\check{\Lambda}=\ZZ^{2n}/\{\check{\epsilon}_i+
\check{\epsilon}_{n+i}-\check{\epsilon}_j-\check{\epsilon}_{n+j}, \; i<j\} 
$$
Let $e_i$ denote the standard basis of $\ZZ^{2n}$. 
The coroots are 
$$
R=\{\pm\alpha_{ij} \; (i< j \in 1,\ldots,n),\;
\pm\beta_{ij}\; (i\le j \in 1,\ldots,n)\},
$$
where $\beta_{ij}=e_i+e_j-e_{n+i}-e_{n+j}$ for
$i<j$ and $\beta_{ii}=e_i-e_{n+i}$. Besides, $\alpha_{ij}=e_i+e_{n+j}-e_j-e_{n+i}$.   

 Fix positive roots 
$$\check{R}^+=\{\check{\alpha}_{ij}
\;\; (i<j\in 1,\ldots,n), \; \check{\beta}_{ij} \;\; (i\le j\in 1,\ldots,n)\}
$$
Then the simple roots are $\check{\alpha}_{12},\ldots,\check{\alpha}_{n-1,n}$ and
$\check{\beta}_{n,n}$.

 For $1\le i< n$ pick the fundamental weight $\check{\omega}_i$
corresponding to the simple coroot $\alpha_{i,i+1}$ to be
$\check{\omega}_i=(1,\ldots,1,0\ldots,0)$, where 1 appears
$i$ times, and 0 appears $2n-i$ times. Let the fundamental weight $\check{\omega}_n$
corresponding to $\beta_{n,n}$ be $\check{\omega}_n=(1,\ldots,1,0,\ldots,0)$, where 1
appears $n$ times. The orthogonal to the coroot latice is the subgroup
$\ZZ\check{\omega}_0$ with $\check{\omega}_0=(1,0,\ldots,0;1,0,\ldots,0)$. 

 All our conditions are satisfied and $\gamma=(1,\ldots,1;0,\ldots,0)$ (here 1 appears $n$
times).

 For $1\le i<n$ let $\gamma_i=(2,\ldots,2,1,\ldots,1;0,\ldots,0,1,\ldots,1)$, where 2 appears
$i$ times then 1 appears $n-i$ times then 0 appears $i$ times and finally 1 appears $n-i$
times. The element $\omega=(1,\ldots,1)\in\Lambda$ generates the group of coweights orthogonal
to all roots. The semigroup $\Lambda^+_{G,S}$ is the $\ZZ_+$-span of
$\gamma,\gamma_1,\ldots,\gamma_{n-1},\omega$. In fact, these $n+1$ elements are linearly
independent in $\Lambda$, so $\Lambda^+_{G,S}\,\iso\,(\ZZ_+)^{n+1}$.

 Note that $V^{\gamma}$ is the spinor representation of $\check{G}\,\iso\,\GSpin_{2n+1}$ of
dimension $2^n$. We have
$$
V^{\gamma}\otimes V^{\gamma}\iso V^{2\gamma}\oplus V^{\omega}\oplus \sum_{i=1}^{n-1}
V^{\gamma_i}
$$
and $\dim V^{\gamma_1}=2n+1$. Besides, $\wedge^i V^{\gamma_1}\,\iso\, V^{\gamma_i+(i-1)\omega}$
for
$i=1,\ldots, n-1$ and $\wedge^n V^{\gamma_1}\,\iso\, V^{2\gamma+(n-1)\omega}$. There is an
exact sequence $1\to \Gm\to \check{G}\to \SO_{2n+1}\to 1$, and $V^{\gamma_1-\omega}$
comes from the standard representation of $\SO_{2n+1}$.

\medskip\noindent
3. {\bf The case $G=\GSpin_{2n+1}$, $n\ge 1$.} The group $G$ is the quotient of
$\Gm\times\Spin_{2n+1}$ by the diagonally embedded $\{\pm 1\}$. We have
$\check{G}\;\iso\,\GSp_{2n}$, the root data for $G$ is dual to that of example~2. Interchanging
the role of objects and coobjects in example 2, we get 
$$
\Lambda=\ZZ^{2n}/\{\epsilon_i+
\epsilon_{n+i}-\epsilon_j-\epsilon_{n+j}, \; i<j\} 
$$
where $(\epsilon_i)$ is the standard basis of $\ZZ^{2n}$. The weight lattice is
$$
\check{\Lambda}=\{(a_1,\ldots,a_{2n})\in\ZZ^{2n}\mid
a_i+a_{n+i} \mbox{ does not depend on } i\}
$$
Define
$\check{\omega}_i$ as follows.
For $1\le i< n$ let
$$
\check{\omega}_i=(2,\ldots,2,1,\ldots,1;0,\ldots,0,1,\ldots,1)\in\check{\Lambda}
$$ 
where 2
appears $i$ times then 1 appears $n-i$ times then 0 appears $i$ times and finally 1 appears
$n-i$ times. Let
$\check{\omega}_n=(1,\ldots,1;0,\ldots,0)$ where 1 appears $n$ times. Let
$\check{\omega}_0=(1,\ldots,1)$.

 All the conditions are satisfied and $\gamma=(1,0,\ldots,0)\in\Lambda$.

For $1\le i\le n$ let $\gamma_i=(1,\ldots,1,0,\ldots,0)\in\Lambda$, where 1 appears $i$ times
and 0 appears $2n-i$ times. The group of coweights orthogonal to all roots is $\ZZ\omega$
with $\omega=(1,0,\ldots,0;1,0,\ldots,0)$, here 1 appears on the first and $(n+1)$-th places.
The semigroup $\Lambda^+_{G,S}$ is the $\ZZ_+$-span of $\gamma_1,\ldots,\gamma_n,\omega$.
These $n+1$ elements are linearly independent in $\Lambda$, so
$\Lambda^+_{G,S}\,\iso\,(\ZZ_+)^{n+1}$. 

 We have $\wedge^2 V^{\gamma}=V^{\gamma_2}\oplus V^{\omega}$. To assure  condition (A) of
Conjecture~\ref{Con_Langlands}, it suffices to require that 
$V^{\gamma}_E$ is irreducible and $V^{\gamma_2}_E$ has no local
subsystems of rank one.

\bigskip\noindent
4. {\bf The case $H=\Spin_{2n}$, $n\ge 3$ odd}. We have
$$
\check{\Lambda}_H=\{(a_1,\ldots,a_n)\in (\frac{1}{2}\ZZ)^n\mid a_i-a_j\in\ZZ\}
$$
and $\Lambda_H=\{(a_1,\ldots,a_n)\in\ZZ^n\mid a_1+\ldots+a_n \;\mbox{is even}\}$. The roots
are
$$
\check{R}=\{\pm\check{\alpha}_{ij}, \pm\check{\beta}_{ij} \; (1\le i< j\le n)\}
$$
with $\check{\alpha}_{ij}=\check{e}_i-\check{e}_j$ and
$\check{\beta}_{ij}=\check{e}_i+\check{e}_j$. The coroots are 
$$
R=\{\pm{\alpha}_{ij}, \pm{\beta}_{ij} \; (1\le i< j\le n)\}
$$
with ${\alpha}_{ij}=e_i-e_j$ and ${\beta}_{ij}=e_i+e_j$. 
Pick positive roots $\check{R}=\{\check{\alpha}_{ij}, \,\check{\beta}_{ij}\; \,
(1\le i< j\le n)\}$. Then simple roots are
$$
\check{\alpha}_{12},\ldots,\check{\alpha}_{n-1,n},\,\check{\beta}_{n-1,n}
$$
We have
$\check{Q}_H=\{(a_1,\ldots,a_n)\in\ZZ^n\mid a_1+\ldots+a_n \;\mbox{is even}\}$,
$$
Q_H=\{(a_1,\ldots,a_n)\in (\frac{1}{2}\ZZ)^n\mid a_i-a_j\in\ZZ\}
$$
There are two possible choices for $\gamma_H$, namely $\gamma_H=(\frac{1}{2},\ldots,\frac{1}{2})$
or $\gamma_H=(\frac{1}{2},\ldots,\frac{1}{2}, -\frac{1}{2})$, the fundamental coweights corresponding
to
$\check{\beta}_{n-1,n}$ or to $\check{\alpha}_{n-1,n}$ respectively.

\bigskip\noindent
5. {\bf The case $E_6$.} So, $H$ is the simply-connected group corresponding to $E_6$ root
system. There are two possible choices for $\gamma_H$, namely $\omega_1$ or $\omega_6$, the
fundamental coweights corresponding to the simple roots $\check{\alpha}_1$ and
$\check{\alpha}_6$ in the Bourbaki table (\cite{Bo}, Ch. 6, no. 4.12, p.219).

\bigskip\noindent
6. {\bf The case $E_7$.} So, $H$ is the simply-connected group corresponding to $E_7$ root
system. Take $\gamma_H$ to be the fundamental coweight $\omega_7$ corresponding to the root 
$\check{\alpha}_7$ in the Bourbaki table (\cite{Bo}, Ch. 6, no. 4.11, p. 217).

\end{document}